\documentclass[lettersize,twocolumn,journal]{IEEEtran}  % Comment this line out
\usepackage{graphicx} % for pdf, bitmapped graphics files
\usepackage{amsmath,bm,times} % assumes amsmath package installed
\usepackage{amssymb} % assumes amsmath package installed
\usepackage{tikz}
\usetikzlibrary{shapes,arrows,backgrounds,fit,positioning}
\usepackage{subfigure}
\usepackage{cite}
\usepackage{soul}
\usepackage{balance}
\usepackage{tabularx}
\usepackage{url}
\allowdisplaybreaks
\usepackage{multirow}
\usepackage{epstopdf}

\newtheorem{propos}{Proposition}
\newtheorem{remark}{Remark}

\newtheorem{lemma}{Lemma}

\usepackage{algorithm}
\usepackage{algpseudocode}

\newcommand{\qed}{$\hfill\blacksquare$}

\title{Semi-Explicit Solution of Some Discrete-Time Mean-Field-Type Games with Higher-Order Costs}

\author{Julian Barreiro-Gomez, \IEEEmembership{Senior Member, IEEE}, Tyrone E. Duncan  \IEEEmembership{Life Fellow Member, IEEE}, \\Bozenna Pasik-Duncan \IEEEmembership{Life Fellow Member, IEEE}, and Hamidou Tembine, \IEEEmembership{Senior Member, IEEE}
\thanks{J. Barreiro-Gomez is with KU Center for Autonomous Robotic Systems, Department of Computer and Information Engineering, Khalifa University, Abu Dhabi 127788, UAE.(e-mail: {\tt\scriptsize julian.barreirogomez@ku.ac.ae}).}
\thanks{H. Tembine is  with Department of Electrical Engineering and Computer Science, School of Engineering  at Universit\'{e} du Qu\'{e}bec \`{a} Trois-Rivi\`{e}res UQTR, Quebec, Canada, and with Timadie, Guinaga, Grabal,  AI Mali,   TF,  WETE, MFTG, LnG Lab, CI4SI;  (e-mail: {\tt\scriptsize tembine@ieee.org}).}
\thanks{T. E. Duncan and B. Pasik-Duncan are with Department of Mathematics, University of Kansas, Lawrence, KS 66044, USA, (e-mail: {\tt\scriptsize duncan@ku.edu, bozenna@ku.edu})}
\thanks{Authors gratefully acknowledge support from TIMADIE, Guinaga, Grabal, LnG Lab for the MFTG for Machine Intelligence project.
Authors gratefully acknowledge support from U.S.
Air Force Office of Scientific Research under grants
number FA9550-17-1-0259.}
\thanks{J. Barreiro-Gomez is profoundly grateful to God and to Our Lady of Lourdes for the blessings of health and life, without which his contributions to this work would not have been possible}
}

\IEEEoverridecommandlockouts

\usepackage{algorithm}
\usepackage{algpseudocode}
\usepackage[mathscr]{euscript}

\begin{document}

\maketitle

\begin{abstract}
Traditional solvable game theory and mean-field-type game theory (risk-aware games) predominantly focus on quadratic costs due to their analytical tractability. Nevertheless, they often fail to capture critical non-linearities inherent in real-world systems. In this work, we present a unified framework for solving discrete-time game problems with higher-order state and strategy costs involving power-law terms. We derive semi-explicit expressions for equilibrium strategies, cost-to-go functions, and recursive coefficient dynamics across deterministic, stochastic, and multi-agent system settings by convex-completion techniques. The contributions include variance-aware solutions under additive and multiplicative noise, extensions to mean-field-type-dependent dynamics, and conditions that ensure the positivity of recursive coefficients. Our results provide a foundational methodology for analyzing non-linear multi-agent systems under higher-order penalization, bridging classical game theory and mean-field-type game theory with modern computational tools for engineering applications.
\end{abstract}


%%
%\newpage 
%\tableofcontents
%\newpage 
\section{Introduction}
Dynamic game theory with quadratic costs has been quite studied and used due to their analytical convenience and well-established solution techniques. As a limitation of this approach, we observe that, nowadays, there are several real-world engineering problems that exhibit non-linearities that go beyond the simple quadratic formulation. Then, the consideration of higher-order costs, e.g., quartic or general power-law functions, allows modeling the aforementioned non-linearities such as state saturations, resource/energy constraints, and risk-sensitivity.

Mean-field-type game (MFTG) theory investigates the interactions of multiple agents including diverse entities such as individuals, animal populations, networked devices, corporations, nations, genetic systems, and others.  Decision-makers within MFTG frameworks can be classified as atomic, non-atomic, or a hybrid of both. The nature of these interactions spans a spectrum from fully cooperative to entirely selfish, including intermediate forms such as partially altruistic, partially cooperative, selfless, co-opetitive, spiteful, and self-abnegating behaviors, as well as combinations thereof \cite{Coopetitive,berge}. MFTG has become a powerful approach to incorporate risk terms such as variance, quantiles, inverse quantiles, skewness, kurtosis, among other, which are not necessarily linear with respect to the probability distribution of individual/common state or action variables. MFTG has demonstrated applicability in diverse domains, including building evacuation, energy systems, next-generation wireless networks, meta-learning, transportation systems, epidemiology, predictive maintenance, network selection, and blockchain technologies. In contrast to classical mean-field games \cite{Jovanovic_2,Lions,Huang}, MFTG does not intrinsically necessitate a large population of agents. Furthermore, MFTG exhibits distinct equilibrium structures and associated Master Adjoint Systems (MASS), attributable, in part, to the non-linearity of expected payoffs with respect to the individual mean-field (the probability measure of an agent's own state) \cite{MAS_HT}.  Semi-explicit solutions for a broad class of MFTGs with non-quadratic cost functions have been derived, facilitating the validation of numerical methods and deep learning algorithms, including transformer-based architectures \cite{nonlinearMFTG}.

In this work, we use convex completion methods to semi-explicitly derive recursive formulations for equilibrium strategies and equilibrium cost-to-go functions that solve non-quadratic discrete-time game problems. Also, we establish conditions for their well-posedness. Our analysis extends to variance-aware games under stochastic perturbations and risk-aware game-theoretic models where selfish agents interact under higher-order cost penalties.
Our results contribute to the broader understanding of non-quadratic mean-field-type games \cite{refmain0,refmain1,refmain2}, providing a structured approach to handling cost functions of arbitrary even order. This is an extension of the work presented in \cite{Higher-order-control}, which exclusively addressed control problems (one decision-maker case), to the game-theoretical problems (multiple decision-makers case). 

\vspace{-0.3cm}
\subsection{Previous Works}

Let us first discuss on previous works related to the computation of semi-explicit solutions to MFTG with quadratic costs, state-variance, and control-variance considerations.

    %\item \textit{Quadratic costs, state-variance, control-variance:} 
Discrete-time linear-quadratic mean-field-type games (LQ-MFTGs) have been widely studied as they provide a mathematical framework for decision-making in small, medium and large-scale multi-agent systems with strategic interactions. These models, which extend classical linear-quadratic control \cite{ref001,ref002,ref003,ref004,ref005} to include mean-field type interactions \cite{ref006}, are used in applications such as finance, energy storage, and networked systems. Recent research has explored solution methods, stability conditions, and numerical approximations including deep neural networks approximators for discrete-time LQ-MFTGs.
In  \cite{ref04}, a class of discrete-time distributed mean-variance problem is introduced and solved in a semi-explicit way.  This work belongs to the class of discrete-time LQ-MFTGs.  Semi-explicit solutions to continuous time LQ-MFTGs are examined in  \cite{ref01,ref02,ref03}.  In \cite{ref0}, a multi-person mean-field-type game is formulated and solved, which is described by a linear jump-diffusion system of mean-field type and a quadratic cost functional involving the second moments, the square of the expected value of the state, and the control actions of all decision-makers. The authors propose a direct method to solve the game, team, and bargaining problems. As an advantage, the direct method solution approach does not require solving the Bellman-Kolmogorov equations or backward-forward stochastic differential equations of Pontryagin's type.  The work in \cite{ref0b} studies the blockchain cryptographic tokens by means of mean-field-type game theory. It introduces the variance-aware utility function per decision-maker to capture the risk of cryptographic tokens associated with the uncertainties of technology adoption, network security, regulatory legislation, and market volatility. The authors provide a semi-explicit solution to the mean-variance problem in blockchain token economics.
The book on \textit{Mean-Field-Type Games for Engineers} \cite{ref1}, presents both continuous-time and discrete-time formulations of mean-field-type games, including LQ-MFTGs. The book discusses convex optimization techniques and equilibrium computations and provides a MatLab toolbox for these problems \cite{Toolbox}. The work in  \cite{ref2} investigates discrete-time linear-quadratic mean-field-type repeated games under different information structures, including perfect, incomplete, and imperfect information settings.
In the energy domain, \cite{ref3}  introduce a stochastic model predictive control framework based on LQ-MFTGs for microgrid energy storage optimization. This work demonstrates how game-theoretic models can enhance decision-making in power systems. In \cite{waterref} a mean-field-type model predictive control, which is a class of risk-aware control considering within the cost functional not only the mean but also the variance of both the system states and control inputs, is considered. The authors applied the semi-explicit solutions to water   distribution networks. The work in \cite{ref4}  employs adaptive dynamic programming techniques to solve LQ-MFTCs for unknown mean-field stochastic discrete-time systems. The study focuses on real-time learning approaches and optimal control strategies.
In \cite{ref4b}, a profit optimization between electricity producers is formulated and solved. The problem is described by a linear jump-diffusion system of conditional mean-field type where the conditioning is with respect to common noise and a quadratic cost functional involving the second moment, the square of the conditional expectation of the control actions of the producers. They provide semi-explicit solution of the corresponding mean-field-type game problem with common noise. Semi-explicit solutions to LQ-MFTGs with multiple control inputs are also discussed in \cite{ref4c,viswa_2023}.

% \begin{table*}[t!]
%     \caption{Summary of Contributions Across Different Problem Types}
%     \label{tab:contributions}
%     \centering
%     \renewcommand{\arraystretch}{1.2}
%     \begin{tabular}{lccc}%{|p{4cm}|p{2.5cm}|p{2.5cm}|p{2.5cm}|}
%          %\\
%         \hline
        
%          \textbf{Problem Type}  & \textbf{Equilibrium State} & \textbf{Equilibrium Strategy (per agent)} & \textbf{Equilibrium Cost (per agent)} 
%          \\
%         \hline
%         Selfish Agents with Higher-Order Costs  & \(\checkmark\) & \(\checkmark\) & \(\checkmark\) \\
%         %\hline
%         Selfish Agents with Variance-Aware Higher-Order Costs &  \(\checkmark\) & \(\checkmark\) & \(\checkmark\) \\
%         %\hline
%         Risk-Aware Games Under Multiplicative Noise &  \(\checkmark\) & \(\checkmark\) & \(\checkmark\) \\
%         \hline
%     \end{tabular}
% \end{table*}

On the other hand, we have the works devoted to the analysis and computation of 
semi-explicit solutions to MFTG with non-quadratic costs, state, and control higher-order moments.

MFTGs with higher-order cost structures extend classical game theory to multi-agent systems where decision-makers optimize objectives involving non-linear or polynomial cost functions. Finding explicit solutions in these settings is challenging due to the complexity of the resulting coupled systems of partial-integro differential equations (PIDEs) or stochastic integro-differential equations (SIDEs). However, recent advances have introduced semi-explicit and fully explicit methods to solve these problems efficiently. The works in  \cite{ref5b} and \cite{ref5c} provide various solvable examples beyond the classical linear-quadratic game problems. These include quadratic-quadratic games and games with power, logarithmic, sine square, hyperbolic sine square payoffs. Non-linear state dynamics such as log-state, control-dependent regime switching, quadratic state, cotangent state, and hyperbolic cotangent state are considered. They identify equilibrium strategies and equilibrium payoffs in state-and-conditional mean-field type feedback form. They show that the direct method can be used to solve broader classes of non-quadratic mean-field-type games under jump-diffusion-regime switching Gauss-Volterra processes which include fractional Brownian motions and multi-fractional Brownian motions. They provide semi-explicit solutions to the fully cooperative, non-cooperative non-zero-sum, and adversarial game problems (see also  \cite{ref5c}).
Hierarchical mean-field-type games with polynomial cost functions are introduced in    \cite{ref5}. Their study develops closed-loop semi-explicit solutions for master systems with hierarchical structures, allowing for leadership and multi-level optimization. In \cite{ref7} the  problem of designing a collection of terminal payoff of mean-field type by interacting decision-makers to a specified terminal measure is considered. The author solves
in a semi-explicit way a class of distributed planning in mean-field-type games with different objective functionals of higher-order and establishes  some relationships between the proposed framework, optimal
transport theory, and distributed control of dynamical systems.
The work in \cite{ref6} examines a class of mean-field-type games with finite number of decision-makers with state dynamics driven by Rosenblatt processes and higher-order cost. Rosenblatt processes are non-Gaussian, non-Poisson, and non-Markov with long-range dependence. The authors provide, in a semi-explicit way, equilibrium strategies and equilibrium costs in linear state-and-mean-field-type feedback form for all decision-makers. They  show that the equilibrium strategies for state driven by Brownian, multi-fractional Brownian, Gauss-Volterra processes no longer provide equilibrium in presence of Rosenblatt processes.

\vspace{-0.3cm}
\subsection{Contributions}
\textit{Selfish Agents with Higher-Order Costs}
\begin{itemize}
    \item We introduce a multi-agent game framework where each agent optimizes a higher-order cost function.
    \item We develop semi-explicit equilibrium strategies and equilibrium costs for selfish agents interacting through a shared system dynamics model.
    \item We establish recursive conditions for the existence of state-feedback  Nash equilibrium strategies with quartic and higher-order costs.
    \item We extend the approach to general even-order cost functions.
\end{itemize}
%\item 

\textit{Selfish Agents with Variance-Aware Higher-Order Costs}
\begin{itemize}
    \item We extend the game-theoretic model to stochastic environments where agents account for both mean-state deviations and variance-aware penalties $var(x)+  \bar{x}^{2p}$,  $var(u_{ik})+  \bar{u}_{ik}^{2p}$ per agent, $x$ being a common state to all agents and $u_{ik}$ being the control action picked by agent $i$ at time step $k$, and $\bar{x}_k= \mathbb{E}[x_k]$ and $\bar{u}_{ik}= \mathbb{E}[u_{ik}]$.
    \item We derive recursive state-and-mean-field-type feedback Nash equilibrium strategies that incorporate both risk sensitivity and higher-order cost minimization.
        \item We establish conditions for the well-posedness of the variance-aware game and present  validation of equilibrium strategies.
\end{itemize}
%\item 
%
\textit{Risk-Aware Games Under Multiplicative Noise}
\begin{itemize}
    \item We develop a stochastic game framework where agents operate in environments with multiplicative noise that is  state-control action-and-mean-field-type dependent. 
    \item We establish semi-explicit state-and-mean-field-type feedback Nash equilibrium strategies and equilibrium costs for risk-sensitive agents with higher-order cost penalties $(x-\bar{x})^{2o}+ \bar{x}^{2p}$,  $(u_{ik}-\bar{u}_{ik})^{2o}+ \bar{u}_{ik}^{2p}$
    $o,p\geq 1.$
    \item We extend the recursive solution framework to general even-order cost functions, ensuring applicability to a broad class of nonlinear game settings.
    \item We provide a verification method using a measure-space formulation, confirming the robustness of the proposed equilibrium solutions.
\end{itemize}
%\end{enumerate}
%

\vspace{-0.3cm}
\subsection{How does the Convex Completion Works?}
As it was presented before in the literature review, there are many works on the computation of semi-explicit solutions to mean-field-type game problems by following the so-called direct method. One of the key elements in this method, in the framework of linear-quadratic settings, is the square completion in order to optimize the Hamiltonian over the control actions. In this work, we do not longer work with the square completion as the proposed costs in this paper go beyond the quadratic structure and examines higher-order costs. Instead, this work uses convex completion as briefly introduced next.

Let $X$ be a real vector space and $X^*$ be its dual space. Let $f: X \to \mathbb{R} \cup \{+\infty\}$ be a proper, convex, and lower semi-continuous function. The Legendre-Fenchel conjugate\footnote{Also known as convex conjugate.} (or transform) of $f$ is defined as:
$$f^*(y) = \sup_{x \in X} \{\langle y, x \rangle - f(x)\},$$
where $\langle y, x \rangle$ denotes the dual pairing between $y \in X^*$ and $x \in X$.
The convexity inequality, also known as Fenchel's inequality or Young's inequality, states that for any $x \in X$ and $y \in X^*$:
$$f(x) + f^*(y) \geq \langle y, x \rangle.$$
The convex completion method uses this inequality as follows. After making the difference between the candidate cost and game cost functional of the agent, one gets a function of the control of that agent, say $f_i(u_i).$ We complete this cost by adding  and removing $f^*_i(v_i) - \langle v_i, u_i \rangle.$
$$ f_i(u_i)+f^*_i(v_i) - \langle v_i, u_i \rangle \geq 0 $$ for control actions $(u_i,v_i).$
The equality holds if and only if $v_i \in \partial f_i(u_i)$, where $\partial f_i(u_i)$ is the sub-differential of $f_i$ at $u_i$. The sub-differential is a generalization of the gradient for non-differentiable functions. 
This extra term added to the difference cost goes to the coefficient matching terms making the solution semi-explicit.

\vspace{-0.3cm}
\subsection{Convexity of Higher-Order Terms}

This paper is going to consider higher order terms, i.e., beyond the quadratic case, with a particular structure. Next, Lemma \ref{lema:highorderterms} shows a convexity property that the development of the paper relies on.
\begin{lemma}
\label{lema:highorderterms}
Let $p\geq 1, a\neq 0, b\neq 0.$ Then the mapping $z \mapsto f(z) = z^{2p} + (az+b)^{2p}$ is strictly convex. \hfill $\square$
\end{lemma}

\vspace{-0.3cm}
\subsection{Structure of the Paper}
The rest of the article is organized as follows. Selfish agents with higher-order costs  are presented in $\S$  \ref{sec:label:4}. In
     $\S$ \ref{sec:label:5}, we present selfish agents with variance-aware higher-order costs.
      $\S$ \ref{sec:label:6} discusses risk-aware games under multiplicative noise. The proofs to all the theoretical results are developed in $\S$ \ref{sec:proofs}. Numerical illustrative examples are presented in $\S$ \ref{sec:examples}. Finally, concluding remarks are drawn in $\S$ \ref{sec:conclusions}.
%
%\subsection{Content of the article}
%

 \section{Selfish Agents with Higher Order Costs}  \label{sec:label:4}
In this section we provide  semi-explicit solution of discrete-time game problems with non-quadratic cost functions, specifically focusing on the minimization of a fourth-order state and control cost of each agent. By using convex-completion techniques, we derive semi-explicit expressions for the equilibrium strategies, the associated equilibrium cost-to-go function, and a recursive relation for the cost coefficients. Furthermore, we establish conditions for the existence of positive coefficients in the cost recursion and extend the method to generalized higher-order costs of the power-law form of  arbitrary even order. The findings provide a foundation for analyzing discrete time games with nonquadratic cost.

\subsection{Selfish Agents with Quartic Costs}
We consider the following discrete-time dynamical system influenced by the decision of $I\geq 2$ agents:
\begin{align}
\label{eq:dynamics_game_1}
    \bar{x}_{k+1} = \bar{a}_k \bar{x}_k +  \sum_{i\in \mathcal{I}}\bar{b}_{ik} \bar{u}_{ik}, 
\end{align}
with $k \in \{0, \dots, N-1\},$ where $\mathcal{I}= \{1, \dots, I\},$ is the set of agents, $\bar{x}_k \in \mathbb{R}$ is the state, $\bar{u}_{ik} \in \mathbb{R}$ is the control input of agent $i$, $\bar{a}_k, \bar{b}_k \in \mathbb{R} $ are system parameters at time step $k.$   The initial state $\bar{x}_0$ is given. The goal is to minimize the cost:
\begin{equation}
\label{eq:game_cost}
    L_i(\bar{x}_0,\bar{u}) = \bar{q}_{iN} \bar{x}_N^4 + \sum_{k=0}^{N-1} \left(\bar{q}_{ik} \bar{x}_k^4 + \bar{r}_{ik} \bar{u}_{ik}^4\right),
\end{equation}
where $\bar{q}_{ik}, \bar{q}_{iN}, \bar{r}_{ik} > 0$ are weights penalizing the state and control deviations. The game problem is given by
\begin{align}
\label{eq:problem5}
    \forall i \in \mathcal{I}:~\min_{\bar{u}_i \in \mathcal{\bar U}_{ik}} \mathbb{E} [L_i(\bar{x}_0,\bar{u})],~\mathrm{s.t.}~\eqref{eq:dynamics_game_1},~\bar{x}_0~\mathrm{given}.
\end{align}
The set of strategies for the $i$-th player at time step $k$ is 
$$\bar{U}_{ik}= \{ \bar{u}_{ik}\in \mathbb{R} | \bar{u}_{ik}^4 < +\infty \}=\mathbb{R},$$  
and the set 
\begin{align}
\label{eq:cal_bar_U_i}
\mathcal{\bar U}_{ik}= \hspace{-0.1cm} \left\{  (\bar{u}_{ik})_{k\in \{0, \dots, N-1\}} \ \hspace{-0.1cm} \bigg| \hspace{-0.1cm} \begin{array}{l}
   \bar{u}_{ik}(.) \in \bar{U}_{ik} , \\
    \bar{u}_{ik}~  \mbox{Lebesgue measurable} 
\end{array} \hspace{-0.15cm} \right\}.    
\end{align}
\begin{propos}
\label{propos:proposition_8}
    The equilibrium strategy for the game problem in \eqref{eq:problem5} with dynamics in \eqref{eq:dynamics_game_1} and cost in \eqref{eq:game_cost} is
    \begin{align}
    \label{eq:optimal_strategy}
        % \bar{u}_{ik}^* = -\frac{\left(\frac{\bar{\alpha}_{i,k+1} \bar{b}_{ik} }{\bar{r}_{ik} }\right)^{1/3}}{1 + \left(\frac{\bar{\alpha}_{i,k+1} \bar{b}_{ik} }{\bar{r}_{ik}}\right)^{1/3} \bar{b}_{ik} }  \left( \bar{a}_k \bar{x}_k +  \sum_{j\in \mathcal{I} \backslash \{ i\}}\bar{b}_{jk} \bar{u}_{jk} \right).\\
        %
        \bar{u}_{ik}^* = - ( \bar{E}^{ -1}_{k} \bar{c}_{k})_i  (\bar{a}_k \bar{x}_k),
    \end{align}
    where $\bar{\alpha}_{ik}$ solves the following recursive equation:
    \begin{align}
    \label{eq:recursive_game}
    \bar{\alpha}_{ik} &= \bar{q}_{ik} + \bar{r}_{ik} \left(( \bar{E}^{ -1} \bar{c})_i  \bar{a}_k\right)^4 \notag\\
    &\quad\quad\quad\quad\quad + \bar{\alpha}_{k+1} \left(\bar{a}_k -  \sum_{j\in \mathcal{I}} ( \bar{E}^{ -1} \bar{c})_j  \bar{a}_k\bar{b}_{jk} 
     \right)^4, 
      \end{align}
      where $\bar{E}_k$ is the square matrix with $ \bar{e}_{ii,k}=1$,  $\bar{e}_{ij,k}=\bar{c}_{ik} \bar{b}_{jk}$,
\begin{align*}
\bar{c}_{ik}= \frac{\left(\frac{\bar{\alpha}_{i,k+1} \bar{b}_{ik} }{\bar{r}_{ik} }\right)^{1/3}}{1 + \left(\frac{\bar{\alpha}_{i,k+1} \bar{b}_{ik} }{\bar{r}_{ik}}\right)^{1/3} \bar{b}_{ik} },
\end{align*}
and $\bar{\alpha}_{iN}= \bar{q}_{iN}$ as terminal condition. \hfill $\square$
\end{propos}

\begin{remark}
Assume that $ \bar{E}$ is invertible, the coefficients $ \bar{\alpha}_{ik}$ remain positive if: \textit{(i)} $\bar{q}_{ik}$, $\bar{q}_{iN}$, $\bar{r}_{ik} > 0$, \textit{(ii)} $\bar{a}_k$, $\bar{b}_{ik}$ are bounded, and \textit{(iii)} denominator in $\bar{u}_{ik}^*$ satisfies that  $1 + \left({{\bar{\alpha}_{i,k+1}} \bar{b}_{ik}}/{\bar{r}_{ik}}\right)^{1/3} \bar{b}_{ik} > 0$.
These conditions are trivially satisfied in most physical systems, ensuring $\bar{\alpha}_{ik} > 0$ for all  $k \in \{0, 1, \ldots, N \}$. \hfill $\square$
\end{remark}

\subsection{Extension to $ \bar{x}^{2p}$ Costs}

We now provide for discrete-time games with individual costs of the form $\bar{x}_k^{2p}$ and $\bar{u}_{ik}^{2p}$, where  $p \geq 1. $ The solutions include the equilibrium control law, the recursive relation for the cost-to-go coefficient, and the expression for the equilibrium cost.
We consider the same discrete-time system as in \eqref{eq:dynamics_game_1}. 
% \begin{equation}
%     \bar{x}_{k+1} = \bar{a}_k \bar{x}_k + \sum_{i\in \mathcal{I}}\bar{b}_{ik}  \bar{u}_{ik} ,
% \end{equation}
% where $\bar{x}_k \in \mathbb{R}$ is the state, $\bar{u}_{ik}  \in \mathbb{R}$ is the control input, and $\bar{a}_k, \bar{b}_{ik}  \in \mathbb{R}$ are system parameters. 
%
The objective is to minimize the $2p$-order cost function of agent $i$:
\begin{align}
\label{eq:cost_game_2}
    L_i(\bar{x}_0,u) = \bar{q}_{iN}  \bar{x}_N^{2p} + \sum_{k=0}^{N-1} \left(\bar{q}_{ik}  \bar{x}_k^{2p} + \bar{r}_{ik}  \bar{u}_{ik} ^{2p}\right),
\end{align}
where $\bar{q}_{ik} , \bar{q}_{iN} , \bar{r}_{ik}  > 0$  are weights. The new set of feasible strategies for the $i$-th player at time step $k$ is 
$\bar{U}_{ik}= \{ \bar{u}_{ik}\in \mathbb{R} | \bar{u}_{ik}^{2p} < +\infty \}=\mathbb{R},$ 
with $\mathcal{\bar U}_k$ as in \eqref{eq:cal_bar_U_i}.
\begin{propos}
\label{propos:proposition_9}
    The equilibrium strategy for the problem in \eqref{eq:problem5} with dynamics in \eqref{eq:dynamics_game_1} and cost functional in \eqref{eq:cost_game_2} is
    \begin{align}
    \label{eq:optimal_strategy_2}
        % &\bar{u}_{ik}^* = \\
        % &-\frac{\left(\frac{{\bar{\alpha}_{i,k+1}} \bar{b}_{ik}}{\bar{r}_{ik}}\right)^{1/(2p-1)} }{1 + \left(\frac{{\bar{\alpha}_{i,k+1}} \bar{b}_{ik}}{\bar{r}_{ik}}\right)^{1/(2p-1)} \bar{b}_{ik}} \left( \bar{a}_k \bar{x}_k + \sum_{j\in \mathcal{I} \backslash \{ i\}}\bar{b}_{jk} \bar{u}_{jk} \right),\notag
        %
        %\bar{u}_{ik}^* &= -\bar{c}_{ik} \left( \bar{a}_k \bar{x}_k + \sum_{j\in \mathcal{I} \backslash \{ i\}}\bar{b}_{jk} \bar{u}_{jk} \right),
        %
        \bar{u}_{ik}^* = - ( \bar{E}^{ -1}_k \bar{c}_k)_i  (\bar{a}_k \bar{x}_k),
    \end{align}
    where $\bar{\alpha}_{i,k}$ solves the following recursive equation:
    \begin{align}
    \label{eq:recursive_game_2}
     \bar{\alpha}_{i,k} &= \bar{q}_{ik} + \bar{r}_{ik} \left( ( \bar{E}^{ -1}_k \bar{c}_k)_i  \bar{a}_k\right)^{2p} \notag\\
     &\quad\quad\quad + \bar{\alpha}_{i,k+1} \left(\bar{a}_k - 
   \sum_{j\in \mathcal{I}} ( \bar{E}^{ -1}_k \bar{c}_k)_j  \bar{a}_k\bar{b}_{jk} \right)^{2p},
    \end{align}
    with $\bar{E}_k$ being a square matrix with $ \bar{e}_{ii}=1$ ,  
$ \bar{e}_{ij}=  \bar{c}_{ik} \bar{b}_{jk}$ and 
\begin{align*}
\bar{c}_{ik} = \frac{\left(\frac{{\bar{\alpha}_{i,k+1}} \bar{b}_{ik}}{\bar{r}_{ik}}\right)^{1/(2p-1)} }{1 + \left(\frac{{\bar{\alpha}_{i,k+1}} \bar{b}_{ik}}{\bar{r}_{ik}}\right)^{1/(2p-1)} \bar{b}_{ik}},    
\end{align*}
and $\bar{\alpha}_{iN}= \bar{q}_{iN}$ as terminal condition. \hfill $\square$
\end{propos}

\section{Selfish Agents with Variance-Aware Higher-Order Costs}  
\label{sec:label:5}

We now focus on variance-aware discrete-time games  with finitely many agents with non-quadratic cost functions, specifically focusing on the minimization of a fourth-order state and control cost per agent. By using convex-completion techniques, we derive semi-explicit expressions for the equilibrium strategies and  the associated variance-aware cost-to-go function, and a recursive relation for the cost coefficients. Furthermore, we establish conditions for the existence of positive coefficients in the variance-aware cost recursion and extend the method to generalized higher-order   costs of the even power-law form. 

\subsection{Selfish Agents with Variance-aware Quartic Costs}
We consider the following discrete-time dynamical system:
\begin{align}
\label{eq:dynamics_game_2}
    {x}_{k+1} = \bar{a}_k {x}_k +   \sum_{i\in \mathcal{I}}\bar{b}_{ik} {u}_{ik} +\epsilon_{k+1} ,
\end{align}
with $k \in \{ 0, \dots, N-1\}$, 
where $\mathcal{I}= \{1, \dots, I\}$ is the set of agents, $\bar{x}_k \in \mathbb{R}$ is the state, $\bar{u}_{ik} \in \mathbb{R}$ is the control input of agent $i$, $\bar{a}_k, \bar{b}_k \in \mathbb{R} $ are system parameters at time step $k.$   The initial state $x_0$ is a given random variable that is independent of $\epsilon.$  $ \{\epsilon_{k}\}_k$ is a zero-mean random process, with $\epsilon_0=0.$ 
The state dynamics can be re-written as 
\begin{align*}
   {x}_{k+1} &=   \bar{a}_k \bar{x}_k +  \sum_{i\in \mathcal{I}}\bar{b}_{ik} \bar{u}_{ik} \\
   &\quad\quad\quad\quad +
   \bar{a}_k (x_k-\bar{x}_k) +  \sum_{i\in \mathcal{I}}\bar{b}_{ik} ({u}_{ik} -\bar{u}_{ik})+\epsilon_{k+1},
\end{align*}
with $k \in \{ 0, \dots, N-1\}$. 
The goal is to minimize the expected value of the cost function:
\begin{align}
\label{eq:cost_game_2}
    L_{i}(x_0,&u) =   {q}_{iN}  {var(x_N)} +\bar{q}_{iN}   \bar{x}_N^4 \\
    &+ \sum_{k=0}^{N-1} \left( {q}_{ik}  {var(x_k)}+\bar{q}_{ik}   \bar{x}_k^4 + {r}_{ik}  {var(u_k)}+ \bar{r}_{ik}   \bar{u}_k^4\right),\notag
\end{align}
where $  q_{ik} ,{q}_{iN} , r_{ik} , \bar{q}_{ik} ,  \bar{q}_{iN},   \bar{r}_{ik}  > 0$ are weights penalizing the state and control deviations. The game problem is:
\begin{align}
\label{eq:problem6}
    \forall i \in \mathcal{I}:~\min_{{u}_i \in \mathcal{U}_{ik}} \mathbb{E} [L_i(x_0,{u})],~\mathrm{s.t.}~\eqref{eq:dynamics_game_2},~x_0~\mathrm{given}.
\end{align}
The set of strategies for the $i$-th player at time step $k$ is 
${U}_{ik}= \{ {u}_{ik}\in \mathbb{R} | \bar{u}_{ik}^4 < +\infty \}=\mathbb{R},$  
and the set 
\begin{align}
\label{eq:cal_U_i}
\mathcal{U}_{ik} \hspace{-0.05cm} = \hspace{-0.05cm} \left\{ \hspace{-0.05cm} ({u}_{ik})_{k\in \{0, \dots, N-1\}} \ \hspace{-0.12cm} \bigg| \hspace{-0.12cm} \begin{array}{l}
   {u}_{ik}(.) \in {U}_{ik} , \\
    {u}_{ik}~  \mbox{Lebesgue measurable} 
\end{array} \hspace{-0.16cm} \right\}.    
\end{align}

\begin{propos}
\label{propos:proposition_10}
    The equilibrium strategy solution for problem in \eqref{eq:problem6} with system dynamics as in \eqref{eq:dynamics_game_2} and cost functional in \eqref{eq:cost_game_2} is given by
    \begin{align}
    \label{eq:optimal_strategy_3}
    u_{i,k}^* = - ( \bar{E}^{ -1}_k \bar{c}_k)_i  (\bar{a}_k \bar{x}_k) -  (E^{-1}_k c_k)_i\bar{a}_k(x_k-\bar{x}_k),
    \end{align}
    with $E_k$ being the matrix with $e_{ii,k}=1$ and $e_{ij,k}= c_{ik}\bar{b}_{jk}$ and $$c_{ik}=\frac{ {\alpha}_{i,k+1} \bar{b}_{ik} }
{r_{ik} +\alpha_{i,k+1} \bar{b}^2_{ik}},$$ 
matrix $\bar{E}_k$ is square with $ \bar{e}_{ii,k}=1$ ,  
$ \bar{e}_{ij,k}=  \bar{c}_{ik} \bar{b}_{jk}$ and $$\bar{c}_{ik}= \frac{\left(\frac{\bar{\alpha}_{i,k+1} \bar{b}_{ik} }{\bar{r}_{ik} }\right)^{1/3}}{1 + \left(\frac{\bar{\alpha}_{i,k+1} \bar{b}_{ik} }{\bar{r}_{ik}}\right)^{1/3} \bar{b}_{ik} },$$
where $\alpha_{ik}$ and $\bar{\alpha}_{ik}$ solve the following recursive equations:
\begin{subequations}
    \label{eq:recursive_game_3}
\begin{align}
    \bar{\alpha}_{ik} &= \bar{q}_{ik} + \bar{r}_{ik} \left(( \bar{E}^{ -1}_k \bar{c}_k)_i  \bar{a}_k\right)^4 \notag\\
    &\quad\quad\quad\, + \bar{\alpha}_{k+1} \left(\bar{a}_k -  \sum_{j\in \mathcal{I}} ( \bar{E}^{ -1}_k \bar{c}_k)_j  \bar{a}_k\bar{b}_{jk} 
     \right)^4, \\
      %\bar{\alpha}_{iN}= \bar{q}_{iN}, \\
  {\alpha}_{ik} &= {q}_{ik} + {r}_{ik} \left(( {E}^{ -1}_k {c}_k)_i  \bar{a}_k\right)^2 \notag\\
  &\quad\quad\quad\, + {\alpha}_{i,k+1} \left(\bar{a}_k -  \sum_{j\in \mathcal{I}} ( {E}^{ -1}_k {c}_k)_j  \bar{a}_k\bar{b}_{jk} 
     \right)^2, %\\
      %{\alpha}_{iN}= {q}_{iN}
\end{align}
\end{subequations}
with $\bar{\alpha}_{iN}= \bar{q}_{iN}$ and ${\alpha}_{iN}= {q}_{iN}$ as terminal conditions. \hfill $\square$
\end{propos}

\begin{remark}
If ${E}_k$ and $\bar{E}_k$ are invertible then  the coefficients ${\alpha}_{ik} ,  \bar{\alpha}_{ik} $ remain positive if: \textit{(i)} ${q}_{ik}$, ${q}_{iN}$, ${r}_{ik}$, $\bar{q}_{ik}$, $\bar{q}_{iN}$, $\bar{r}_{ik}> 0$, \textit{(ii)} $\bar{a}_k$, $\bar{b}_{ik}$ are bounded, \textit{(iii)} denominator in $ \bar{u}_{ik} ^*$ satisfies that  $1 + \left({\bar{\alpha}_{i,k+1} \bar{b}_{ik} }/{\bar{r}_{ik} }\right)^{1/3} \bar{b}_{ik}  > 0$, and \textit{(iv)} denominator in $ u_{ik} ^*-\bar{u}_{ik} ^*$ satisfies that  $r_{ik}  +\alpha_{i,k+1} \bar{b}^2_{ik} >0$. Then $\alpha_{ik}  > 0$  and  $ \bar{\alpha}_{ik} $ for all  $k \in \{0, 1, \ldots, N \}.$ \hfill $\square$
\end{remark}
\subsection{Extension to $ \bar{x}^{2p}$ Costs}

We now address discrete-time games with costs of the form $x_{k}^{2p}$ and $u_{ik}^{2p}$, where  $p \geq 1. $ The solutions include the equilibrium strategy law, the recursive relation for the cost-to-go coefficient, and the expression for the equilibrium cost.
We consider the same discrete-time system as in \eqref{eq:dynamics_game_2}. 
% \begin{equation}
%    {x}_{k+1} = \bar{a}_k {x}_k + \ \sum_{i\in \mathcal{I}}\bar{b}_{ik}  {u}_{ik}+\epsilon_{k+1}, \quad k \in \{ 0, \dots, N-1\},
% \end{equation}
% where $x_k \in \mathbb{R}$ is the state, $u_{ik}  \in \mathbb{R}$ is the control input, and $a_k, \bar{b}_{ik}  \in \mathbb{R}$ are system parameters. 
%
The objective is to minimize the variance-aware $2p$-order cost:
\begin{align}
\label{eq:cost_game_2p}
    L_i&(x_0,u) =  {q}_{iN} {var(x_N)} +
    \bar{q}_{iN} \bar{x}_N^{2p} \\
    &+ \sum_{k=0}^{N-1} \left(  {q}_{ik} {var(x_k)} +
    \bar{q}_{ik} \bar{x}_k^{2p} + {r}_{ik} {var(u_{ik})}+
     \bar{r}_{ik} \bar{u}_{ik}^{2p}\right),\notag
\end{align}
where ${q}_{ik}, {q}_{iN}, {r}_{ik} , \bar{q}_{ik}, \bar{q}_{iN}, \bar{r}_{ik} > 0$  are weights. The new set of strategies for the $i$-th player at time step $k$ is 
${U}_{ik}= \{ {u}_{ik}\in \mathbb{R} | \bar{u}_{ik}^{2p} < +\infty \}=\mathbb{R},$ with $\mathcal{U}_{ik}$ as in \eqref{eq:cal_U_i}.

\begin{propos}
\label{propos:proposition_11}
    The equilibrium strategy solution for problem in \eqref{eq:problem6} with system dynamics as in \eqref{eq:dynamics_game_2} and new cost functional in \eqref{eq:cost_game_2p} is given by
    \begin{align}
    \label{eq:optimal_strategy_4}
        u_{i,k}^*  = - ( \bar{E}^{ -1}_k \bar{c}_k)_i  (\bar{a}_k \bar{x}_k) - (E^{-1}_k c_k)_i \bar{a}_k(x_k-\bar{x}_k)
    \end{align}
    with $E_k$ being a matrix with $e_{ii,k}=1,  e_{ij,k}= c_{ik}\bar{b}_{jk}$, 
    \begin{align*}
        c_{ik}=\frac{\alpha_{i,k+1} \bar{b}_{ik}}{r_{ik} +\alpha_{i,k+1} \bar{b}^2_{ik}};
    \end{align*} 
    and  $\bar{E}_k$ is the matrix with $ \bar{e}_{ii,k}=1$ ,  
$ \bar{e}_{ij,k}=  \bar{c}_{ik} \bar{b}_{jk}$ and 
\begin{align*}
\bar{c}_{ik}= \frac{\left(\frac{\bar{\alpha}_{i,k+1} \bar{b}_{ik}}{\bar{r}_{ik}}\right)^{1/(2p-1)} }{1 + \left(\frac{\bar{\alpha}_{i,k+1} \bar{b}_{ik}}{\bar{r}_{ik}}\right)^{1/(2p-1)} \bar{b}_{ik}},
\end{align*}
with the following recursive equations:
\begin{subequations}
    \label{eq:recursive_4}
\begin{align} 
   \bar{\alpha}_{i,k} &= \bar{q}_{ik} + \bar{r}_{ik} \left( ( \bar{E}^{ -1}_k \bar{c}_k)_i  \bar{a}_k\right)^{2p} \notag\\
   &+ \bar{\alpha}_{i,k+1} \left(\bar{a}_k - 
   \sum_{j\in \mathcal{I}} ( \bar{E}^{ -1} \bar{c})_j  \bar{a}_k\bar{b}_{jk} \right)^{2p},
     \\
      %\bar{\alpha}_{iN}= \bar{q}_{iN}\\
       %\bar{\gamma}_{i0}&= \sum_{k=0}^{N-1}  {{\alpha}_{i,k+1}}  \mathbb{E}\left( \epsilon_{k+1}^{2}  \right), \\ 
       \bar{\gamma}_{i,k} &= \bar{\gamma}_{i,k+1} +  {\alpha}_{i,k+1} \mathbb{E}\left(  \epsilon_{k+1}^2  \right), \\
       %\bar{\gamma}_{iN}&= 0, \\
        {\alpha}_{ik} &= {q}_{ik} + {r}_{ik}
         \left(
         ( {E}^{ -1}_k {c}_k)_i  \bar{a}_k
         \right)^2 \notag\\
         &+ {\alpha}_{k+1}
          \left(\bar{a}_{k}-  \sum_{j}( {E}^{ -1}_k {c}_k)_j  \bar{a}_k \bar{b}_{jk}
          \right)^2, %\\
      %{\alpha}_{iN}&= {q}_{iN}. %, \\
\end{align}
\end{subequations}
and with $\bar{\alpha}_{iN}= \bar{q}_{iN}$, $\bar{\gamma}_{iN}= 0$, and ${\alpha}_{iN}= {q}_{iN}$. \hfill $\square$
\end{propos}

\section{Risk-Aware Games under Multiplicative Noise}  \label{sec:label:6}
Non-quadratic costs, such as $ \bar{x}^4, (x-\bar{x})^4, \bar{u}^4$ or $ (x-\bar{u})^4$, introduce another class of  realistic representation of these non-linearities, emphasizing large deviations in state and control.
We now explicitly solve the variance-aware discrete-time games with higher-order costs. We derive semi-explicit solutions for the  variance-aware equilibrium, the   variance-aware equilibrium cost, and the recursive coefficients for the  variance-aware cost-to-go function. Additionally, we extend the methodology to address general higher-order costs $(x_k-\bar{x}_k)^{2o}$ for arbitrary $o \geq 1.$

\subsection{Variance-aware Quartic Cost under Multiplicative Noise}
We consider the following discrete-time dynamical system with multiplicative state-and-mean state dependent noise:
\begin{align}
    \label{eq:dynamics_game_3}
    {x}_{k+1} = \bar{a}_k {x}_k +   \sum_{i\in \mathcal{I}}\bar{b}_{ik} {u}_{ik}
      + (x_k-\bar{x}_k)\epsilon_{k+1} , 
\end{align}
with $k \in \{ 0, \dots, N-1\}$, where $\mathcal{I}= \{1, \dots, I\},$ is the set of agents, $\bar{x}_k \in \mathbb{R}$ is the state, $\bar{u}_{ik} \in \mathbb{R}$ is the control input of agent $i$, $\bar{a}_k, \bar{b}_k \in \mathbb{R} $ are system parameters at time step $k.$   The initial state $x_0$ is a given random variable that is independent of $\epsilon.$  $ \{\epsilon_{k}\}_k$ is a zero-mean random process, with $\epsilon_0=0.$ 
The state dynamics can be re-written as 
\begin{align*}
   {x}_{k+1} &=   \bar{a}_k \bar{x}_k +  \sum_{i\in \mathcal{I}}\bar{b}_{ik}  \bar{u}_{ik} +
   \bar{a}_k (x_k-\bar{x}_k) \notag\\
   &+  \sum_{i\in \mathcal{I}} \bar{b}_{ik}  (u_{ik} -\bar{u}_{ik} )+ (x_k-\bar{x}_k)\epsilon_{k+1},
\end{align*}
with $k \in \{ 0, \dots, N-1\}$. 
The goal is to minimize the expected value of the cost function:
\begin{align}
\label{eq:cost_game_3}
    L_i&(x_0,u) =   {q}_{iN}  {var(x_N)} +\bar{q}_{iN}   \bar{x}_N^4 \\
    &+ \sum_{k=0}^{N-1} \left( {q}_{ik}  {var(x_k)}+\bar{q}_{ik}   \bar{x}_k^4 + {r}_{ik}  {var(u_{ik})}+ {\bar{r}_{ik}}  \bar{u}_{ik} ^4\right),\notag
\end{align}
where $  q_{ik} ,{q}_{iN} , r_{ik} , \bar{q}_{ik} ,  \bar{q}_{iN} ,   \bar{r}_{ik}  > 0$ are weights penalizing the state and control deviations. The game problem is:
\begin{align}
\label{eq:problem7}
    \forall i \in \mathcal{I}:~\min_{{u}_i \in \mathcal{U}_{ik}} \mathbb{E} [L_i(x_0,{u})],~\mathrm{s.t.}~\eqref{eq:dynamics_game_3},~x_0~\mathrm{given}.
\end{align}
The set of strategies for the $i$-th player at time step $k$ is 
${U}_{ik}= \{ {u}_{ik}\in \mathbb{R} | \bar{u}_{ik}^{4} < +\infty \}=\mathbb{R},$ with $\mathcal{U}_{ik}$ as in \eqref{eq:cal_U_i}.

\begin{propos}
\label{propos:proposition_12}
    The equilibrium strategy for the problem in \eqref{eq:problem7} is given by 
    \begin{align}
    \label{eq:optimal_strategy_5}
        {u}_{ik}^* =  -  (E^{-1}_k c_k)_i \bar{a}_k(x_k-\bar{x}_k)  -   ( \bar{E}^{ -1}_k \bar{c}_k)_i  \bar{a}_k \bar{x}_k,
    \end{align}
    where $E_k$ the matrix with $e_{ii,k}=1$, $e_{ij,k}= c_{ik}\bar{b}_{jk}$ and 
    \begin{align*}
     c_{ik}=\frac{ {\alpha}_{i,k+1} \bar{b}_{ik}  }
{r_{ik} +\alpha_{i,k+1} \bar{b}^2_{ik}};   
    \end{align*}
    $ \bar{E}_k$ is the matrix with $ \bar{e}_{ii,k}=1$ ,  
$ \bar{e}_{ij,k}=  \bar{c}_{ik} \bar{b}_{jk}$ and 
\begin{align*}
\bar{c}_{ik}= \frac{\left(\frac{\bar{\alpha}_{i,k+1} \bar{b}_{ik} }{\bar{r}_{ik} }\right)^{1/3}}{1 + \left(\frac{\bar{\alpha}_{i,k+1} \bar{b}_{ik} }{\bar{r}_{ik}}\right)^{1/3} \bar{b}_{ik} } .    
\end{align*}
The terms $\alpha_{ik}$, $\bar{\alpha}_{ik}$ solve the following recursive equations: 
\begin{subequations}
\label{eq:recursive_5}
\begin{align} 
    \bar{\alpha}_{ik} &= \bar{q}_{ik} + \bar{r}_{ik} \left( ( \bar{E}^{ -1}_k \bar{c}_k)_i  \bar{a}_k\right)^4 \\
    &+ \bar{\alpha}_{i,k+1} \left(\bar{a}_{k}-     \sum_{j\in \mathcal{I} }( \bar{E}^{ -1}_k \bar{c}_k)_j  \bar{a}_k \bar{b}_{jk} \right)^4, \notag\\
      %\bar{\alpha}_{iN} = \bar{q}_{iN} , \\
      {\alpha}_{ik} &= {q}_{ik} + {r}_{ik} \left(( {E}^{ -1}_k {c}_k)_i  \bar{a}_k \right)^2 \\
    &+ {\alpha}_{i,k+1} \left(
    \bar{a}_{k}-     \sum_{j\in \mathcal{I} }( {E}^{ -1}_k {c}_k)_j  \bar{a}_k \bar{b}_{jk} 
    \right)^2+  {\alpha}_{i,k+1} \mathbb{E}\left( \epsilon_{k+1}^2  \right), \notag%\\
    %{\alpha}_{iN} = {q}_{iN}.
\end{align}
\end{subequations}
with terminal conditions $\bar{\alpha}_{iN} = \bar{q}_{iN}$ and ${\alpha}_{iN} = {q}_{iN}$. \hfill $\square$
\end{propos}

\begin{remark}
%\subsubsection{Existence of Positive $\alpha_{ik}$ and $ \bar{\alpha}_{ik}$}
If $E_k$ and $\bar E_k$ are invertible, then the coefficients ${\alpha}_k,  \bar{\alpha}_k$ remain positive if: \textit{(i)} ${q}_{ik}$, ${q}_{iN}$, ${r}_{ik}$, $\bar{q}_{ik}$, $\bar{q}_{iN}$, $\bar{r}_{ik} > 0$, \textit{(ii)} $ \bar{a}_k$, $\bar{b}_{ik}$ are bounded, \textit{(iii)} denominator in $ \bar{u}_{ik}^*$ satisfies  $1 + \left({\bar{\alpha}_{i,k+1} \bar{b}_{ik}}/{\bar{r}_{ik}}\right)^{1/3} \bar{b}_{ik} > 0$, and \textit{(iv)} denominator in $ u_{ik}^*-\bar{u}_{ik}^*$ satisfies  $r_k +\alpha_{i,k+1} \bar{b}^2_{ik}>0$. Then $\alpha_{ik} > 0$  and  $ \bar{\alpha}_{ik}>0$ for all  $k \in \{0, 1, \ldots, N \}$. \hfill $\square$
\end{remark}
\subsection{Extension to $ \bar{x}^{2p}$ Costs}
We now provide for discrete-time  game problems with costs of the form $\bar x_k^{2p}$ and $\bar u_{ik}^{2p}$, where  $p \geq 1. $ The solutions include the equilibrium control law, the recursive relation for the cost-to-go coefficient, and the expression for the equilibrium cost. We consider the discrete-time system in \eqref{eq:dynamics_game_3}. 
% \begin{equation}
%    {x}_{k+1} = \bar{a}_k {x}_k + \sum_{i\in \mathcal{I}}\bar{b}_{ik}  {u}_{ik}+(x_k-\bar{x}_k)\epsilon_{k+1}, \quad k \in \{ 0, \dots, N-1\},
% \end{equation}
% where $x_k \in \mathbb{R}$ is the state, $u_{ik}  \in \mathbb{R}$ is the control input, and $a_k, \bar{b}_{ik}  \in \mathbb{R}$ are system parameters. 
The new objective is to minimize the variance-aware $2p$-order cost:
\begin{align}
\label{eq:cost_game_4}
    L_i&(x_0,u)=  {q}_{iN} {var(x_N)} +
    \bar{q}_{iN} \bar{x}_N^{2p} \\
    &+ \sum_{k=0}^{N-1} \left(  {q}_{ik} {var(x_k)} +
    \bar{q}_{ik} \bar{x}_k^{2p} + {r}_{ik} {var(u_{ik})}+
     \bar{r}_{ik} \bar{u}_{ik}^{2p}\right), \notag
\end{align}
where ${q}_{ik}, {q}_{iN}, {r}_{ik} , \bar{q}_{ik}, \bar{q}_{iN}, \bar{r}_{ik} > 0$  are weights. The set of strategies for the $i$-th player at time step $k$ is 
${U}_{ik}= \{ {u}_{ik}\in \mathbb{R} | \bar{u}_{ik}^{2p} < +\infty \}=\mathbb{R},$ with $\mathcal{U}_{ik}$ as in \eqref{eq:cal_U_i}.
\begin{propos}
\label{propos:proposition_13}
    The equilibrium strategy for the problem in \eqref{eq:problem7} with the cost functional in \eqref{eq:cost_game_4} is
    \begin{align}
        \label{eq:optimal_strategy_6}
        {u}_{ik}^* =   - ( {E}^{ -1}_k {c}_k)_i  \bar{a}_k (x_k-\bar{x}_k)- ( \bar{E}^{ -1}_k \bar{c}_k)_i  \bar{a}_k \bar{x}_k
    \end{align}
where $E_k$ is a matrix with $e_{ii,k}=1,  e_{ij,k}= c_{ik}\bar{b}_{jk}$, 
\begin{align*}
    c_{ik}=\frac{\alpha_{i,k+1} \bar{b}_{ik}}
{r_{ik} +\alpha_{i,k+1} \bar{b}^2_{ik}},
\end{align*}
and $\bar{E}_k$ is the square matrix with $\bar{e}_{ii}=1$, $\bar{e}_{ij,k}=  \bar{c}_{ik} \bar{b}_{jk}$ and 
\begin{align*}
\bar{c}_{ik}= -\frac{\left(\frac{\bar{\alpha}_{i,k+1} \bar{b}_{ik}}{\bar{r}_{ik}}\right)^{1/(2p-1)} }{1 + \left(\frac{\bar{\alpha}_{i,k+1} \bar{b}_{ik}}{\bar{r}_{ik}}\right)^{1/(2p-1)} \bar{b}_{ik}}.   
\end{align*}
The terms $\bar{\alpha}_{ik}$, and ${\alpha}_{ik}$ solve the recursive equations
\begin{subequations}
\label{eq:recursive_6}
    \begin{align}
        \bar{\alpha}_{ik} &= \bar{q}_{ik} + \bar{r}_{ik} \left( ( \bar{E}^{ -1}_k \bar{c}_k)_i  \bar{a}_k\right)^{2p} \\
    &+ \bar{\alpha}_{k+1} \left(\bar{a}_{k}-     \sum_{j\in \mathcal{I} }( \bar{E}^{ -1}_k \bar{c}_k)_j  \bar{a}_k \bar{b}_{jk} \right)^{2p}, \notag\\
      %\bar{\alpha}_{iN} &= \bar{q}_{iN} , \\
      {\alpha}_{ik} &= {q}_{ik} + {r}_{ik} \left(( {E}^{ -1}_k {c}_k)_i  \bar{a}_k \right)^2 \\
    &+ {\alpha}_{i,k+1} \left(
    \bar{a}_{k}-     \sum_{j\in \mathcal{I} }( {E}^{ -1}_k {c}_k)_j  \bar{a}_k \bar{b}_{jk} 
    \right)^2+  {\alpha}_{i,k+1} \mathbb{E}\left( \epsilon_{k+1}^2  \right). \notag%\\
    %{\alpha}_{iN} &= {q}_{iN},
    \end{align}
\end{subequations}
with terminal conditions $\bar{\alpha}_{iN} = \bar{q}_{iN}$ and ${\alpha}_{iN} = {q}_{iN}$. \hfill $\square$
\end{propos}

\subsection{Extension to $ \mathbb{E}[(x-\mathbb{E}{x})^{2o}]$ Costs with  Multiplicative state-control and mean-field type dependent noise}
We now provide for discrete-time games  with $2o$-moment costs of the form $\mathbb{E}[(x_k-\bar{x}_k)^{2o}]$ and $\mathbb{E}[(u_{ik}-\bar{u}_{ik})^{2o}]$, where  $o \geq 1. $ The solutions include the equilibrium control law, the recursive relation for the cost-to-go coefficient, and the expression for the equilibrium costs. 
%
%\subsubsection{ State dynamics with Multiplicative state-control and mean-field dependent noise}
We consider the following discrete-time system of mean-field type:
\begin{align}
\label{eq:dynamics_game_4}
   {x}_{k+1} &=  \left(\bar{a}_k \bar{x}_k + \sum_{i\in \mathcal{I}}\bar{b}_{ik} \bar{u}_{ik} \right) \notag\\
   &+
   \left({a}_k ({x}_k- \bar{x}_k) + \sum_{i\in \mathcal{I}} {b}_{ik}  ({u}_{ik} -\bar{u}_{ik} )\right)\epsilon_{k+1},
\end{align}
with $k \in \{ 0, \dots, N-1\}$, 
where $x_k \in \mathbb{R}$ is the state, $u_{ik}  \in \mathbb{R}$ is the control input, and $\bar{a}_k, \bar{b}_{ik}  \in \mathbb{R}$ are system parameters. 
%
%
%\subsubsection{ Cost functional of  agent $i$ }
The objective is to minimize the  $2o$-th moment-dependent  cost:
\begin{align}
\label{eq:cost_game_5}
    L_i(x_0,u) &=  {q}_{iN}  {\mathbb{E}[(x_N-\bar{x}_N)^{2o}]}+
    \bar{q}_{iN} \bar{x}_N^{2p} \notag\\
    &+ \sum_{k=0}^{N-1}  {q}_{ik}  {\mathbb{E}[(x_k-\bar{x}_k)^{2o}]} +
    \bar{q}_{ik} \bar{x}_k^{2p}\notag\\
    &+  \sum_{k=0}^{N-1} 
     {r}_{ik}  {\mathbb{E}[(u_{ik} -\bar{u}_{ik})^{2o}]}+
     \bar{r}_{ik}  \bar{u}_{ik}^{2p},
\end{align}
where ${q}_{ik}, {q}_{iN} , {r}_{ik} , \bar{q}_{ik}, \bar{q}_{iN}, \bar{r}_{ik}  > 0$  are weights. The game theory problem is stated as follows:
\begin{align}
\label{eq:problem8}
    \forall i \in \mathcal{I}:~\min_{{u}_i \in \mathcal{U}_{ik}} \mathbb{E} [L_i(x_0,{u})],~\mathrm{s.t.}~\eqref{eq:dynamics_game_4},~x_0~\mathrm{given}.
\end{align}
The set of strategies for the $i$-th player at time step $k$ is 
${U}_{ik}= \{ {u}_{ik}\in \mathbb{R} | \bar{u}_{ik}^{2p} < +\infty \}=\mathbb{R},$ with $\mathcal{U}_{ik}$ as in \eqref{eq:cal_U_i}.
\begin{propos}
\label{propos:proposition_14}
    The equilibrium strategy for the problem in \eqref{eq:problem8} with the cost functional in \eqref{eq:cost_game_5} is
    \begin{align}
        \label{eq:optimal_strategy_7}
        u_{ik}^*= -    (\tilde E^{-1}_k \tilde c_k)_i a_k({x}_k-\bar{x}_k) -   (\bar E^{-1}_k \bar c_k)_i\bar{a}_k \bar{x}_k,
    \end{align}
    where $\tilde{E}_k$ is a matrix with $\tilde{e}_{ii,k}=1$, $\tilde{e}_{ij,k}= \tilde c_{ik} b_{jk}$, and
    \begin{align*}
        \tilde c_{ik}=  \frac{ \left(\frac{{\alpha}_{i,k+1}b_ {ik}m_{k+1, 2o} }{r_{ik}}\right)^{{1}/{(2o-1)}}}{1+ \left(\frac{{\alpha}_{i,k+1}b_{ik}m_{k+1, 2o} }{r_{ik}}\right)^{{1}/{(2o-1)}}b_{ik}},
    \end{align*}
and $\bar{E}_k$ is the square matrix with $\bar{e}_{ii}=1$, $\bar{e}_{ij,k}=  \bar{c}_{ik} \bar{b}_{jk}$ and 
\begin{align*}
\bar{c}_{ik}= -\frac{\left(\frac{\bar{\alpha}_{i,k+1} \bar{b}_{ik}}{\bar{r}_{ik}}\right)^{1/(2p-1)} }{1 + \left(\frac{\bar{\alpha}_{i,k+1} \bar{b}_{ik}}{\bar{r}_{ik}}\right)^{1/(2p-1)} \bar{b}_{ik}}.   
\end{align*}
The terms ${\alpha}_{i,k}$ and $\bar{\alpha}_{i,k}$ solve the recursive equations
\begin{subequations}
\label{eq:recursive_7}
    \begin{align}
        \bar{\alpha}_{ik} &= \bar{q}_{ik} + \bar{r}_{ik} \left( ( \bar{E}^{ -1}_k \bar{c}_k)_{i}  \bar{a}_k\right)^{2p} \notag\\
    &+ \bar{\alpha}_{k+1} \left(\bar{a}_{k}-     \sum_{j\in \mathcal{I} }( \bar{E}^{ -1}_k \bar{c}_k)_j  \bar{a}_k \bar{b}_{jk} \right)^{2p}, \\
      %\bar{\alpha}_{iN} = \bar{q}_{iN} , \\
      {\alpha}_{ik} &= {q}_{ik} + {r}_{ik} \left(( \tilde{E}^{ -1}_k \tilde{c}_k)_i  {a}_k \right)^{2o} \notag\\
    &+ {\alpha}_{i,k+1} \left(
    {a}_{k}-     \sum_{j\in \mathcal{I} }( \tilde{E}^{ -1}_k \tilde{c}_k)_j  {a}_k {b}_{jk} 
    \right)^{2o} , %\\
    %{\alpha}_{iN} = {q}_{iN}.
    \end{align}
\end{subequations}
with $\bar{\alpha}_{iN} = \bar{q}_{iN}$ and ${\alpha}_{iN} = {q}_{iN}$; and $m_{k+1,2o}  = \mathbb{E}[\epsilon_{k+1}^{2o}]$. \hfill $\square$
\end{propos}

\subsection{Verification Method in the Space of Measures}
The above analysis covers noise beyond the Markovian noise. 
We now verify the proposed solution by checking the case where the noise process is Markovian. To do so we work directly on the space of probability measures associated with $x_k$ and denote by $\mu_k= \mathbb{P}_{x_k} $ its probability distribution. Then $\mu_{k+1}$ can be obtained  $\mathbb{P}_{\epsilon_{k+1} },$   $\mu_{k}$ and $x_k, \bar{x}_k, u_{ik}, \bar{u}_{ik}$. Note here that the mean-field type terms are involved in the transition kernel.  The cost can also be written in terms of   $\mu_{k}$ and the functional $u_{ik},$ of mean-field type.
Our new augmented state becomes $\mu$  and the system becomes a big deterministic system in the space of measures. As an advantage, one can write a dynamic programming principle of this new augmented state. 
\begin{align*}
i &\in \mathcal{I}=\{1, \ldots, I \}, \\
V_{ik}&(\mu_k)= \inf_{u_{ik}} \Bigg(  {q}_{ik}  {\int \left[(z-\int y \mu_k(dy))^{2o} \right] \mu_k(dz) } \\
&+
    \bar{q}_{ik} \left(\int y \mu_k(dy)\right)^{2p}  \\ 
    &+ {r}_{ik} {\int \left[ \left(u_{ik}(z,\mu_k)-\int  u_{ik}(y, \mu_k) \mu_k(dy) \right)^{2o} \right] \mu_k(dz) \Bigg) }\\
    &+
     \bar{r}_{ik}  \left(\int u_{ik}(y,\mu_k) \mu_k(dy)\right)^{2p}  + V_{i,k+1}(\mu_{k+1}) \Bigg),
\end{align*}
and the terminal cost of agent $i$  is 
\begin{align*}
    V_{i,N}(\mu_N)&={q}_{iN} {\int \left[\left(z-\int y \mu_N(dy)\right)^{2o}\right] \mu_N(dz)}\\
    &+
    \bar{q}_{iN} \left(\int y \mu_N(dy)\right)^{2p},\\
%\end{align*}
%where 
%\begin{eqnarray} \nonumber
\mu_{k+1}(dy) &= \left\{ \hspace{-0.2cm}  
\begin{array}{ll}
  \mathbb{P}_{\epsilon_{k+1} }\bigg( \left({a}_k ({x}_k-\bar{x}_k) + \sum\limits_{i\in \mathcal{I} }{b}_{ik} ({u}_{ik}-\bar{u}_{ik})\right)^{-1}\\
  \left(dy-  \left(\bar{a}_k \bar{x}_k + \sum\limits_{i\in \mathcal{I} } \bar{b}_{ik} \bar{u}_{ik}\right) \right) \bigg) \ \\ 
  \mbox{if } \left({a}_k ({x}_k-\bar{x}_k) + \sum\limits_{i\in \mathcal{I} }{b}_{ik} ({u}_{ik}-\bar{u}_{ik}) \right)\neq 0,\\ \\
\delta_{\bar{a}_k \bar{x}_k + \sum_{i\in \mathcal{I} }\bar{b}_{ik} \bar{u}_{ik}} (dy) \ \mbox{otherwise },
\end{array} 
\right.
\end{align*}
We seek for $V_{ik}(\mu_k)$ in the form  
\begin{align*}
\alpha_{ik}  {\int \left[\left(z-\int y \mu_k(dy)\right)^{2o}\right] \mu_k(dz)}+  \bar{\alpha}_{ik} \left(\int y \mu_k(dy)\right)^{2p}.   
\end{align*}
Inserting the above yields the same dynamics for ${\alpha}_{ik}$ and $\bar{\alpha}_{ik}$ as well as the equilibrium  structure
\begin{align*}
u_{ik}^*(x, \mu_k) &= -    (\tilde E^{-1}_k \tilde c_k)_i a_k\left({x}-\int y \mu_k(dy)\right) \\
&-   (\bar E^{-1}_k \bar c_k)_i\bar{a}_k \left(\int y \mu_k(dy)\right).    
\end{align*}
%%%%%%%%
In the next Section, we present all the proofs corresponding to the aforementioned results in Lemma \ref{lema:highorderterms} and Proposition \ref{propos:proposition_8}-Proposition \ref{propos:proposition_14}.

\section{Proofs}
\label{sec:proofs}

\begin{proof}[Proof of Lemma \ref{lema:highorderterms}]
The function $f$ is twice differentiable as it is polynomial. It is strictly convex if its second derivative $f''(z) > 0$, for all $z$.
We differentiate twice to obtain:
\[
f''(z) 
= 2p(2p-1)z^{2p-2} + 2p(2p-1)a^2 (az+b)^{2p-2}
\]
We observe:
\begin{itemize}
\item When $p\geq 1,$ the term $2p(2p-1)z^{2p-2}$ is \textit{positive} for all $z \neq 0$.
\item When $p\geq 1, a, b \neq 0$ The term $2p(2p-1)a^2 (az+b)^{2p-2}$ is always \textit{positive}, for $z \neq - b/a \neq 0$. 
\end{itemize}
Since both terms are non-negative and do not vanish simultaneously , the function is strictly convex under the assumption that $p\geq 1, a, b \neq 0.$ This complete the proof. \qed
\end{proof}

\vspace{0.1cm}
\begin{proof}[Proof of Proposition \ref{propos:proposition_8}]
    We propose a candidate cost  function at time $k$ to be:
\begin{align*}
    f_{ik}(\bar{x}_k) = \bar{{\alpha}}_{ik} \bar{x}_k^4,
\end{align*}
where $\bar{{\alpha}}_{ik} > 0$ is a scalar coefficient to be determined recursively. Using the telescopic sum yields
\begin{align*}
    f_{i,N}(\bar{x}_N) = f_{i,0}(\bar{x}_0)+ \sum_{k=0}^{N-1}  f_{i,k+1}(\bar{x}_{k+1})-f_{i,k}(\bar{x}_{k}).
\end{align*}
By taking the difference $L_i(x_0,u) - f_{i,0}(\bar{x}_0)$ we obtain: 
\begin{align*}
 L_i&(x_0,u) - f_{i,0}(\bar{x}_0) = \bar{q}_{iN}\bar{x}_N^4-f_{iN}(\bar{x}_N) \\
 &+ \sum_{k=0}^{N-1} \left(\bar{q}_{ik} \bar{x}_k^4 + \bar{r}_{ik} \bar{u}_{ik}^4\right)+ f_{i,k+1}(\bar{x}_{k+1})-f_{i,k}(\bar{x}_{k}).
\end{align*}
Introducing the term:
\begin{equation*}
    C_{ik} = \min_{\bar{u}_{ik}} \left( \bar{r}_{ik} \bar{u}_k^4 + {\bar{\alpha}_{i,k+1}} \left(\bar{a}_k \bar{x}_k + \sum_{i\in \mathcal{I}}\bar{b}_{ik} \bar{u}_{ik}\right)^4 \right),
\end{equation*}
and replacing back in the difference equation yields
\begin{align*}
   L_i(x_0,u) &- f_{i,0}(\bar{x}_0) = \bar{q}_{iN} \bar{x}_N^4-f_{i,N}(\bar{x}_N) \\
   &+ \sum_{k=0}^{N-1} \left( \bar{r}_{ik} \bar{u}_{ik}^4+ f_{i,k+1}(\bar{x}_{k+1})\right)-C_{ik} \\
   &+  \left(C_{ik}+\bar{q}_{ik} \bar{x}_k^4 -f_{i,k}(\bar{x}_{k})\right).
\end{align*}
To find the best response strategy, we differentiate the cost function with respect to $\bar{u}_k$:
\begin{align*}
   4 \bar{r}_{ik} \bar{u}_{ik}^3 + 4 \bar{\alpha}_{i,k+1} \left(\bar{a}_k \bar{x}_k +  \sum_{j\in \mathcal{I}}\bar{b}_{jk} \bar{u}_{jk}\right)^3 \bar{b}_{ik} = 0.
\end{align*}
As the terms are continuously differentiable and convex, the first order condition (Fermat's rule) is also sufficient.
We simplify and isolate  $(\bar{a}_k \bar{x}_k +  \sum_{i\in \mathcal{I}}\bar{b}_{ik} \bar{u}_{ik})^3$, i.e.,
\begin{equation*}
    \bar{u}_{ik}^3 = -\frac{{\bar{\alpha}_{i,k+1}} \bar{b}_{ik} }{\bar{r}_{ik} } \left(\bar{a}_k \bar{x}_k + \sum_{j\in \mathcal{I}}\bar{b}_{jk} \bar{u}_{jk}\right)^3.
\end{equation*}
Taking the cube root one obtains 
\begin{align*}
    \bar{u}_{ik}^* = -\bar{c}_{ik}  \left( \bar{a}_k \bar{x}_k +  \sum_{j\in \mathcal{I} \backslash \{ i\}}\bar{b}_{jk} \bar{u}_{jk} \right),
\end{align*}
which means that $ \bar{E}_k  \bar{u}_k =  -  \bar{c}_k  (\bar{a}_k \bar{x}_k)$. In addition, if the matrix $E_k$ is invertible then it  follows that  the equilibrium strategy of agent $i$ is given by \eqref{eq:optimal_strategy}.
Substitute $\bar{u}_{ik}^*$ back into the difference  equation and collect terms proportional to $\bar{x}_k^4$, and one gets the recursive equation for $\bar{\alpha}_{ik}$, completing the proof. \qed
\end{proof}

\vspace{0.1cm}
\begin{proof}[Proof of Proposition \ref{propos:proposition_9}]
The candidate cost of agent $i$  is assumed to have the following  form:
\begin{align*}
    f_{ik} (\bar{x}_k) = \bar{\alpha}_{ik}  \bar{x}_k^{2p},
\end{align*}
where $\bar{\alpha}_{ik} > 0$ is a coefficient to be determined recursively. The difference term now has
\begin{align*}
 C_{ik}  = \min_{\bar{u}_{ik} } \left( \bar{r}_{ik}  \bar{u}_{ik} ^{2p} + {\bar{\alpha}_{i,k+1}} \left(\bar{a}_k \bar{x}_k + \sum_{j\in \mathcal{I}}\ \bar{b}_{jk}  \bar{u}_{jk} \right)^{2p} \right).
\end{align*}
Differentiating with respect to $\bar{u}_{ik}$ gives:
\begin{align*}
   2p \, \bar{r}_{ik} \bar{u}_{ik}^{2p-1} + 2p \, {\bar{\alpha}_{i,k+1}} \bar{b}_{ik} \left( \bar{a}_k \bar{x}_k + \sum_{j\in \mathcal{I}}\ \bar{b}_{jk}  \bar{u}_{jk}  \right)^{2p-1} = 0.
\end{align*}
Rearrange terms:
\begin{align*}
    \bar{r}_{ik} \bar{u}_{ik}^{2p-1} + {\bar{\alpha}_{i,k+1}} \bar{b}_{ik} \left(\bar{a}_k \bar{x}_k + \sum_{j\in \mathcal{I}}\ \bar{b}_{jk}  \bar{u}_{jk} \right)^{2p-1} = 0.
\end{align*}
We isolate  $\left(\bar{a}_k \bar{x}_k +  \sum_{j\in \mathcal{I}}\ \bar{b}_{jk}  \bar{u}_{jk}\right)^{2p-1}$:
\begin{align*}
    \bar{u}_{ik}^{2p-1} = -\frac{{\bar{\alpha}_{i,k+1}} \bar{b}_{ik}}{\bar{r}_{ik}} \left( \bar{a}_k \bar{x}_k +  \sum_{j\in \mathcal{I}}\ \bar{b}_{jk}  \bar{u}_{jk} \right)^{2p-1}.
\end{align*}
Take the $(2p-1)$-th root on both sides one lands at
\begin{align*}
    \bar{u}_{ik}^* &= -\bar{c}_{ik} \left( \bar{a}_k \bar{x}_k + \sum_{j\in \mathcal{I} \backslash \{ i\}}\bar{b}_{jk} \bar{u}_{jk} \right),
\end{align*}
which means that $ \bar{E}_k  \bar{u}_k =  -  \bar{c}_k  (\bar{a}_k \bar{x}_k)$ where $\bar{E}_k$ is the announced square matrix. If the matrix $E_k$ is invertible then it  follows that  the equilibrium strategy of agent $i$ is as in \eqref{eq:optimal_strategy_2}. 
%\begin{align*}
%    \bar{u}_{ik}^* = - ( \bar{E}^{ -1}_k \bar{c}_k)_i  (\bar{a}_k \bar{x}_k)
%\end{align*}
Substituting the best response strategy $\bar{u}_{ik}^*$ back into the candidate cost yields the recursive equation for $\bar{\alpha}_{i,k}$ in \eqref{eq:recursive_game_2}, %The equilibrium cost at each time step is proportional to $\bar{x}_k^{2p}$:
%\begin{equation*}
    %$f_{ik} (\bar{x}_k) = \bar{\alpha}_{ik}  \bar{x}_k^{2p}$, 
%\end{equation*}
completing the proof. \qed
\end{proof}

\vspace{0.1cm}
\begin{proof}[Proof of Proposition \ref{propos:proposition_10}]
    We propose a candidate cost  function at time $k$ to be:
\begin{align*}
    f_{ik} (x_k, \bar{x}_k) =  {\alpha}_{ik}   {var(x_k)} +\bar{\alpha}_{ik}   \bar{x}_k^4+\bar\gamma_{ik} ,
\end{align*}
where ${\alpha}_{ik}, \bar{\alpha}_{ik}, \bar\gamma_{ik}   \geq   0$ are  scalar coefficients to be determined recursively. Computing by telescopic sum yields
\begin{align*}
    f_{iN}(x_N) = f_{i0}(x_0)+ \sum_{k=0}^{N-1}  f_{i,k+1}(x_{k+1})-f_{i,k}(x_{k})  .
\end{align*}
By taking the difference $L_i(x_0,u) - f_{i0}(x_0)$ we obtain: 
\begin{align*}
   L_i&(x_0,u) - f_{i0}(x_0) =  {q}_{iN} {var(x_N)} +\bar{q}_{iN} \bar{x}_N^4-f_N(x_N)
   \\
   &+ \sum_{k=0}^{N-1} \left( {q}_{ik} {var(x_k)}+\bar{q}_{ik} \bar{x}_k^4 + {r}_{ik} {var(u_k)}+ \bar{r}_{ik} \bar{u}_{ik}^4\right)\\
   &+ f_{i,k+1}(x_{k+1})-f_{i,k}(x_{k}).
\end{align*}
Let us introduce the following term:
\begin{align*}
    C_{ik} &= \min_{u_{ik}} 
    \Bigg( {r}_{ik} {var(u_{ik}) } +  {\alpha}_{i,k+1}  {var(x_{k+1}) } 
    +\bar{r}_{ik} \bar{u}_{ik}^4 \\
    &+ \bar{\alpha}_{i,k+1} (\bar{a}_{k}\bar{x}_k + \sum_{j\in \mathcal{I}} \bar{b}_{jk} \bar{u}_{jk})^4 \Bigg),
\end{align*}
and replacing back in the difference equation one arrives at
\begin{align*}
   L_i&(x_0,u)- f_{i0}(x_0) =  {q}_{iN} {var(x_N)} +
    \bar{q}_{iN} \bar{x}_N^4-f_{iN}(x_N) \\ 
    &+ \sum_{k=0}^{N-1} \left( {r}_{ik} {var(u_{ik})}+\bar{r}_{ik} \bar{u}_{ik}^4+ f_{i,k+1}(x_{k+1})\right)-C_{ik}\\ 
   &+  \left(C_{ik}+ {q}_{ik} {var(x_k)}+\bar{q}_{ik} \bar{x}_k^4 -f_{i,k}(x_{k})\right),
\end{align*}
and we have that
\begin{align*}
&{r}_{ik} {var(u_{ik})} +  {\alpha}_{i,k+1}  {var(x_{k+1})} = \mathbb{E} \Bigg[{r}_{ik} {(u_{ik}- \bar{u}_{ik})^2}  \\
&+  {\alpha}_{i,k+1} 
\Bigg( \bar{a}_k (x_k-\bar{x}_k) +  \sum_{j\in \mathcal{I}}\bar{b}_{jk} (u_{jk}-\bar{u}_{jk})+\epsilon_{k+1} \Bigg)^2 \Bigg],
\end{align*}
 which means that 
% 
%${r}_k { \color{blue}  (u_k- \bar{u}_k) }  = -  {\alpha}_{k+1} \bar{b}_k  ( \bar{a}_k (x_k-\bar{x}_k) + \bar{b}_k (u_k-\bar{u}_k))$
\begin{align*}
(u_{ik}- \bar{u}_{ik})  &= -  c_{ik} \left(\bar{a}_k(x_k-\bar{x}_k) + \hspace{-0.2cm} \sum_{j\in \mathcal{I} \backslash \{ i\}} \hspace{-0.2cm} \bar{b}_{jk} ({u}_{jk}-\bar{u}_{jk})\right), \notag
\end{align*}
and using $E_k$ matrix we have  
$E_k (u_{k}- \bar{u}_{k})= -   c_k[\bar{a}_k(x_k-\bar{x}_k)]$ whenever $E_k$ is invertible.
To find the deterministic part of best response strategy $ \bar{u}_{ik}^*$ of agent $i,$  we differentiate the cost function with respect to $u_{ik}$, i.e.,
\begin{align*}
   4 \bar{r}_{ik} \bar{u}_{ik}^3 + 4 \bar{\alpha}_{i,k+1} \left(\bar{a}_{k}\bar{x}_k + \sum_{j\in \mathcal{I} }\ \bar{b}_{jk} \bar{u}_{jk}\right)^3 \bar{b}_{ik} = 0.
\end{align*}
We simplify and isolate $\left(\bar{a}_{k}\bar{x}_k +  \sum_{j\in \mathcal{I} }\  \bar{b}_{jk} \bar{u}_{jk}\right)^3$:
\begin{align*}
    \bar{u}_{ik}^3 = -\frac{{\bar{\alpha}_{i,k+1}} \bar{b}_{ik} }{\bar{r}_{ik} } \left(\bar{a}_{k}\bar{x}_k + \ \sum_{j\in \mathcal{I} }\  \bar{b}_{jk} \bar{u}_{jk} \right)^3.
\end{align*}
Taking the cube root:
\begin{equation*}
    \bar{u}_{ik}^* = - \bar{c}_{ik} \left( \bar{a}_k \bar{x}_k +
    \sum_{j\in \mathcal{I} \backslash \{ i\}}\bar{b}_{jk} \bar{u}_{jk} \right) .
\end{equation*}
This means that $\bar{E}_k  \bar{u}_k =  -  \bar{c}_k  (\bar{a}_k \bar{x}_k)$ where $ \bar{E}_k$ is the aforementioned square matrix. If the matrix $\bar{E}_k$ is invertible then it  follows that  the equilibrium strategy of agent $i$ is given by 
\begin{equation*}
    \bar{u}_{ik}^* = - ( \bar{E}^{ -1}_k \bar{c}_k)_i  (\bar{a}_k \bar{x}_k).
\end{equation*}
Now, substitute ${u}_{ik}^*-\bar{u}_{ik}^*$, and $\bar{u}_{ik}^*$ back into the difference  equation and collect terms proportional to $(x_k-\bar{x}_k)^4, \bar{x}_k^4$. Thus, one obtains the announced recursive equations for $\alpha_{ik}$ and $\bar{\alpha}_{ik}$, completing the proof. \qed
\end{proof}

\vspace{0.1cm}
\begin{proof}[Proof of Proposition \ref{propos:proposition_11}]
    The candidate cost is assumed to have the form:
\begin{align*}
    f_{ik}(x_k, \bar{x}_k) = {\alpha}_{ik} var({x}_k)+\bar{\alpha}_{ik} \bar{x}_k^{2p} + \bar{\gamma}_{ik},
\end{align*}
where $\alpha_{ik}, \bar{\alpha}_{ik},   \bar{\gamma}_{ik}\geq  0$ is a coefficient to be determined recursively. The expected value of difference term now has
\begin{align*}    
 C_{ik} &= \min_{u_{ik}} \Bigg(   {r}_{ik} {var(u_{ik})} +  {\alpha}_{i,k+1}  {var(x_{k+1})} +
 \bar{r}_{ik} \bar{u}_{ik}^{2p} \\
 &+ \bar{\alpha}_{i,k+1} (\bar{a}_{k}\bar{x}_k + \sum_{j\in \mathcal{I}}\ \bar{b}_{jk} \bar{u}_{jk})^{2p} \Bigg).
\end{align*}
The stochastic part of the best-response is determined by 
\begin{align*}
{r}_{ik} &{var(u_{ik})} +  {\alpha}_{i,k+1}  {var(x_{k+1})} 
= \mathbb{E} ({r}_{ik} {(u_{ik}- \bar{u}_{ik})^2}  \\
&+  {\alpha}_{i,k+1}  
\left( \bar{a}_k (x_k-\bar{x}_k) +  \sum_{j\in \mathcal{I}}\bar{b}_{jk} (u_{jk}-\bar{u}_{jk})+\epsilon_{k+1} )^2\right),
\end{align*}
 which means that 
 %${r}_k { \color{blue}  (u_k- \bar{u}_k) }  = -  {\alpha}_{k+1} \bar{b}_k  ( \bar{a}_k (x_k-\bar{x}_k) + \bar{b}_k (u_k-\bar{u}_k))$
\begin{align*}
 (u_k- \bar{u}_k)  = - E^{-1}_kc_k \bar{a}_k(x_k-\bar{x}_k) 
\end{align*}
 with the announced matrix $E_k$. The deterministic part of the best response strategy  is obtained as follows:
\begin{align*}
    &\bar{u}_{ik} =-\bar{c}_{ik} \left( \bar{a}_k \bar{x}_k + \sum_{j\in \mathcal{I} \backslash \{ i\}}\bar{b}_{jk} \bar{u}_{jk} \right).
\end{align*}
This means that $ \bar{E}_k  \bar{u}_k =  -  \bar{c}_k  (\bar{a}_k \bar{x}_k)$ where $\bar{E}_k$ is the announced square matrix. If the matrix $\bar{E}_k$ is invertible then it  follows that  the equilibrium strategy of agent $i$ is given by 
\begin{equation*}
    \bar{u}_{ik}^* = - ( \bar{E}^{ -1}_k \bar{c}_k)_i  (\bar{a}_k \bar{x}_k),
\end{equation*}
and the announced equilibrium strategy \eqref{eq:optimal_strategy_4} is obtained. Substituting the equilibrium strategy $u_{ik}^*$ back into the candidate cost yields one arrives at the announced recursive equations, where we have used that  
\begin{align*} 
    (X+ \epsilon_{k+1})^{2} = X^{2} +\epsilon_{k+1}^{2} +  2 X \epsilon_{k+1}.
\end{align*}
Finally, the equilibrium cost  of agent $i$ is 
\begin{align*}
        f_{i0}(x_0, \bar{x}_0) = {\alpha}_{i0} var({x}_0) + \bar{\alpha}_{i0} \bar{x}_0^{2p}+ \bar\gamma_{i0},
\end{align*}
completing the proof. \qed
\end{proof}

\vspace{0.1cm}
\begin{proof}[Proof of Proposition \ref{propos:proposition_12}]
We propose a candidate cost function at time $k$ to be:
\begin{align*}
    f_{ik} (x_k, \bar{x}_k) =  {\alpha}_{ik}   {var(x_k)}+\bar{\alpha}_{ik}   \bar{x}_k^4,
\end{align*}
where $ {\alpha}_{ik} , \bar{\alpha}_{ik}   \geq   0$ are  scalar coefficients to be determined recursively. Using the telescopic sum yields
\begin{align*}
    f_{iN} (x_N) = f_{i0} (x_0)+ \sum_{k=0}^{N-1}  f_{i,k+1}(x_{k+1})-f_{i,k}(x_{k}).
\end{align*}
By taking the difference $L_i(x_0,u) - f_{i0} (x_0)$ we obtain: 
\begin{align*} 
   L_i(x_0,u)- f_{i0}&(x_0) =  {q}_{iN}  {var(x_N)} +\bar{q}_{iN}  \bar{x}_N^4-f_{iN} (x_N)
   \\
   &+ \sum_{k=0}^{N-1} \big( {q}_{ik}  {var(x_k) }+\bar{q}_{ik}  \bar{x}_k^4 + {r}_{ik}  {var(u_k)}\\
   &+ \bar{r}_{ik}  \bar{u}_k^4\big)+ f_{i,k+1}(x_{k+1})-f_{i,k}(x_{k}).
\end{align*}
Considering the following term:
\begin{align*}
    C_{ik}  &= \min_{u_{ik} } 
    \Bigg( {r}_{ik}  {var(u_{ik} ) } +  {\alpha}_{i,k+1}  {var(x_{k+1}) } 
    +\bar{r}_{ik}  \bar{u}_{ik} ^4 \\
    &+ \bar{\alpha}_{i,k+1} \left(\bar{a}_{k}\bar{x}_k +  \sum_{j\in \mathcal{I}}\bar{b}_{jk} \bar{u}_{jk} \right)^4 \Bigg),
\end{align*}
the difference equation can be rewritten as
\begin{align*} \nonumber 
   L_i(x_0,u) &- f_{i0} (x_0) =  {q}_{iN}  {var(x_N)} +
    \bar{q}_{iN}  \bar{x}_N^4-f_N(x_N) \\ 
    &+ \sum_{k=0}^{N-1} \left( {r}_{ik}  {var(u_{ik})}+\bar{r}_{ik}  \bar{u}_{ik}^4+ f_{i,k+1}(x_{k+1})\right)\\
    &-C_{ik}  
   +  \left(C_{ik} + {q}_{ik}  {var(x_k)}+\bar{q}_{ik}  \bar{x}_k^4 -f_{i,k}(x_{k})\right).
\end{align*}
We have that
\begin{align*}
{r}_{ik} {var(u_{ik})} &+  {\alpha}_{i,k+1}  {var(x_{k+1})} 
= \mathbb{E} \Bigg[{{r}_{ik} {(u_{ik}- \bar{u}_{ik})^2}}  \\
&+  {{\alpha}_{i,k+1}}  
\Bigg( \bar{a}_k (x_k-\bar{x}_k) + \sum_{j\in \mathcal{I}}\bar{b}_{jk} ({u}_{jk}-\bar{u}_{jk}) \\
&+(x_k-\bar{x}_k)\epsilon_{k+1}  \Bigg)^2\Bigg],
\end{align*}
 which means that 
 \begin{align*}
&{u_{ik}- \bar{u}_{ik}}  = -  c_{ik} \left( \bar{a}_k  (x_k-\bar{x}_k) +   \hspace{-0.2cm} \sum_{j\in \mathcal{I}\backslash \{ i\}}\ \hspace{-0.2cm} \bar{b}_{jk} ({u}_{jk}-\bar{u}_{jk})  \right).     
 \end{align*}
Using the introduced matrix $E_k$, it follows that  $E_k (u_{k}- \bar{u}_{k})= -   c_k\bar{a}_k(x_k-\bar{x}_k)$ whenever $E_k$ is invertible. 
To find the deterministic part of best response strategies $ \bar{u}_{ik}^*$, we differentiate the cost function with respect to $u_{ik}:$
\begin{align*}
   4 \bar{r}_{ik} \bar{u}_{ik}^3 + 4 \bar{\alpha}_{i,k+1} \left( \bar{a}_{k}\bar{x}_k + \sum_{j\in \mathcal{I} }\ \bar{b}_{jk} \bar{u}_{jk} \right)^3 \bar{b}_{ik} = 0.
\end{align*}
We simplify and isolate  $\left(\bar{a}_{k}\bar{x}_k +  \sum_{j\in \mathcal{I} }\  \bar{b}_{jk} \bar{u}_{jk}\right)^3$:
\begin{align*}
    \bar{u}_{ik}^3 = -\frac{{\bar{\alpha}_{i,k+1}} \bar{b}_{ik} }{\bar{r}_{ik} } \left(\bar{a}_{k}\bar{x}_k + \ \sum_{j\in \mathcal{I} }\  \bar{b}_{jk} \bar{u}_{jk} \right)^3.
\end{align*}
Taking the cube root:
\begin{align*}
    \bar{u}_{ik}^* = -
    \bar{c}_{ik} \left(\bar{a}_k \bar{x}_k +
    \sum_{j\in \mathcal{I} \backslash \{ i\}}\bar{b}_{jk} \bar{u}_{jk} \right).
\end{align*}
This means that $\bar{E}_k  \bar{u}_k =  -  \bar{c}_k  (\bar{a}_k \bar{x}_k)$ where $ \bar{E}_k$ is the announced square matrix. If the matrix $\bar{E}_k$ is invertible then it  follows that  the equilibrium strategy of agent $i$ is given by 
\begin{align*}
    \bar{u}_{ik}^* = - ( \bar{E}^{ -1}_k \bar{c}_k)_i  (\bar{a}_k \bar{x}_k).
\end{align*}
Then, the announced equilibrium strategy in \eqref{eq:optimal_strategy_5} of agent $i$ is obtained. Substitute $u_{ik}^*$ back into the difference  equation and collect terms proportional to $(x_k-\bar x_k)^4$, the recursive equations in \eqref{eq:recursive_5} are  obtained, completing the proof. \qed
\end{proof}

\vspace{0.1cm}
\begin{proof}[Proof of Proposition \ref{propos:proposition_13}]
The candidate cost  of agent $i$ is assumed to have the form:
\begin{align*}
    f_{ik}(x_k, \bar{x}_k) = {\alpha}_{ik} var({x}_k)+\bar{\alpha}_{ik} \bar{x}_k^{2p},
\end{align*}
where $\alpha_{ik}, \bar{\alpha}_{ik},   \bar{\gamma}_{ik}\geq  0$ is a coefficient to be determined recursively. The expected value of difference term now has
\begin{align*}
 C_{ik} &= \min_{u_{ik}} \Bigg(   {r}_{ik} {var(u_{ik})} +  {\alpha}_{i,k+1}  {var(x_{k+1})} +
 \bar{r}_{ik} \bar{u}_{ik}^{2p} \\
 &+ \bar{\alpha}_{i,k+1} (\bar{a}_{k}\bar{x}_k + \sum_{j\in \mathcal{I}}\ \bar{b}_{jk} \bar{u}_{jk})^{2p} \Bigg).
\end{align*}
The stochastic part of the best response is determined by 
\begin{align*}
{r}_{ik} {var(u_{ik})} &+  {\alpha}_{i,k+1}  {var(x_{k+1})} 
= \mathbb{E} ({r}_{ik} {(u_{ik}- \bar{u}_{ik})^2}  \\
&+  {\alpha}_{i,k+1}  
\Bigg( \bar{a}_k (x_k-\bar{x}_k) \\
&+  \sum_{j\in \mathcal{I}}\bar{b}_{jk} (u_{jk}-\bar{u}_{jk})+(x_k-\bar{x}_k)\epsilon_{k+1} \Bigg)^2,
\end{align*}
 which means that 
$E_k(u_k- \bar{u}_k)  = - c_k  \bar{a}_k(x_k-\bar{x}_k)$ with $E_k$ being the announced matrix. The deterministic part of the best response strategy of agent $i$ is obtained as follows:
\begin{align*}
    \bar{u}_{ik} = -\bar{c}_{ik} \left(\bar{a}_k \bar{x}_k + \hspace{-0.2cm} \sum_{j\in \mathcal{I} \backslash \{ i\}}\hspace{-0.2cm}\bar{b}_{jk} \bar{u}_{jk} \right).
\end{align*}
This means that $ \bar{E}_k  \bar{u}_k =  -  \bar{c}_k  (\bar{a}_k \bar{x}_k)$ where $\bar{E}_k$ is the announced square matrix. If the matrix $\bar{E}_k$ is invertible then it follows that  the equilibrium strategy of agent $i$ is given by 
\begin{equation*}
    \bar{u}_{ik}^* = - ( \bar{E}^{ -1}_k \bar{c}_k)_i  (\bar{a}_k \bar{x}_k).
\end{equation*}
Combining together the equilibrium  strategy of agent $i$ yields the equilibrium strategy in \eqref{eq:optimal_strategy_6}. Substituting the equilibrium strategy $u_{ik}^*$ back into the candidate cost yields the recursive equations in \eqref{eq:recursive_6}, where we have used that  
\begin{align*} 
    (a Y+ Y\epsilon_{k+1})^{2} = Y^2(a^{2} +\epsilon_{k+1}^{2} +  2 a\epsilon_{k+1}).
\end{align*}
The equilibrium cost of $i$  is 
\begin{align*}
    f_{i0} (x_0, \bar{x}_0) = {\alpha}_{i0}  var({x}_0) + \bar{\alpha}_{i0} \bar{x}_0^{2p}.
\end{align*}
As the variance of the noises $ \mathbb{E}\left( \epsilon_{k+1}^2  \right)$ vanish, we retrieve the solution of  deterministic game. \qed
\end{proof}

\vspace{0.1cm}
\begin{proof}[Proof of Proposition \ref{propos:proposition_14}]
    The candidate cost is assumed to have the form:
\begin{align*}
    f_{ik} (x_k, \bar{x}_k) =  {\alpha}_{ik}  {\mathbb{E}[(x_k-\bar{x}_k)^{2o}]} +\bar{\alpha}_{ik}  (\bar{x}_k)^{2p},
\end{align*}
where $\alpha_{ik},   \bar\alpha_{ik} \geq  0$ is a coefficient to be determined recursively. 
The expected value of difference term now has
\begin{align*}
 C_{ik} &= \min_{u_{ik} } \Bigg(   {r}_{ik}  {\mathbb{E}[(u_{ik} -\bar{u}_{ik})^{2o}] } +  {\alpha}_{i,k+1}  {\mathbb{E}[(x_{k+1}-\bar{x}_{k+1})^{2o}]}\\
 &
 +\bar{r}_{ik}  \bar{u}_{ik}^{2p} + \bar{\alpha}_{i,k+1} \left(\bar{a}_{k}\bar{x}_k +  \sum_{j\in \mathcal{I}}\ \bar{b}_{jk}  \bar{u}_{jk} \right)^{2p} 
   \Bigg).
\end{align*}
This minimization is subject to $\mathcal{F}_k$ adaptivity of equilibrium strategy. We compute 
\begin{align*}
\mathbb{E}&\left[(x_{k+1}-\bar{x}_{k+1})^{2o}\right] \\
&=   
\mathbb{E}\left[ \left( \left({a}_k ({x}_k-\bar{x}_k) +  \sum_{i\in \mathcal{I}}\ {b}_{ik} ({u}_{ik}-\bar{u}_{ik})\right)\epsilon_{k+1} \right)^{2o}\right]  \\
&=  \mathbb{E}\left[\left({a}_k ({x}_k-\bar{x}_k) +  \sum_{i\in \mathcal{I}}\ {b}_{ik} ({u}_{ik}-\bar{u}_{ik})\right)^{2o}\right]  m_{k+1, 2o}.
\end{align*}
The stochastic part of the  best response of agent $i$  yields
\begin{align*}
    u_{ik}-\bar{u}_{ik} &= - \tilde{c}_{ik}
    \left(a_k({x}_k-\bar{x}_k) + \hspace{-0.2cm} \sum_{j\in \mathcal{I} \backslash \{i\}} \hspace{-0.2cm} {b}_{jk} ({u}_{jk}-\bar{u}_{jk})\right).
\end{align*}
\begin{figure}[t!]
    \centering
    \begin{tabular}{cc}
       \includegraphics[width=0.5\columnwidth]{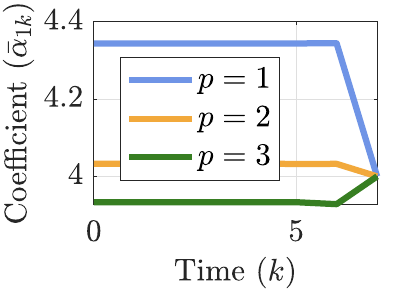}  &
       \hspace{-0.5cm}
       \includegraphics[width=0.5\columnwidth]{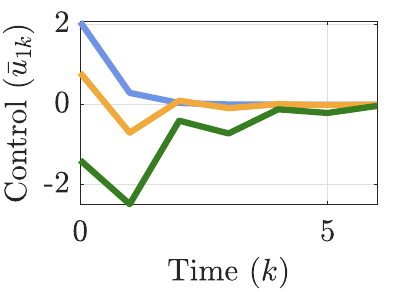}  \\
         \includegraphics[width=0.5\columnwidth]{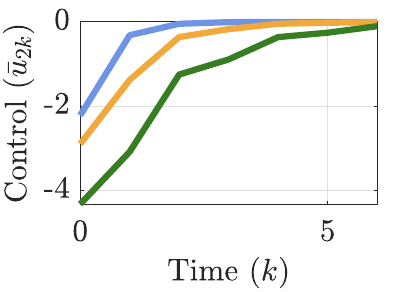}  &
       \hspace{-0.5cm}
       \includegraphics[width=0.5\columnwidth]{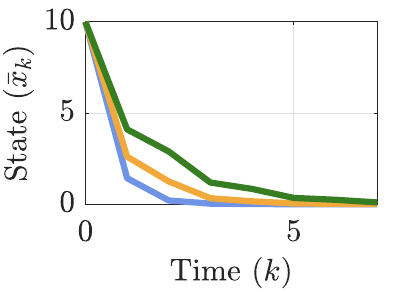}  
    \end{tabular}
    \caption{Results example corresponding to Proposition \ref{propos:proposition_8} and Proposition \ref{propos:proposition_9}.}
    \label{fig:p1_2_games}

    \vspace{-0.3cm}
\end{figure}
Then 
%\begin{align*}
$\tilde{E}_k (u_{k}-\bar{u}_{k}) =  \tilde c_k a_k({x}_k-\bar{x}_k)$,
%\end{align*}
with the announced matrix $\tilde{E}_k$. The deterministic part of the control yields
\begin{align*}
&\bar{u}_{ik}= - \bar{c}_{ik} \left(\bar{a}_k \bar{x}_k + 
\sum_{j\in \mathcal{I} \backslash \{i\}}\bar{b}_{jk} \bar{u}_{jk}
\right).    
\end{align*}
We can easily check that the obtained equilibrium strategy of agent $i$ is the presented in \eqref{eq:optimal_strategy_7}, and it is adapted to the filtration generated by $\{ x_0, \bar{x}_0, \ldots, x_k, \bar{x}_k \}$. Substituting the equilibrium strategy    $u_{ik}^*$ back into the candidate cost of $i$ yields the recursive equations in \eqref{eq:recursive_7}. Note that  the  ${\alpha}_{ik}$ and ${\alpha}_{jk}$  are coupled through $\tilde{E}_k$ and $\tilde{c}_k.$ 
The $\bar{\alpha}_{ik}$ and $\bar{\alpha}_{jk}$  are coupled through $\bar{E}_k$ and $\bar{c}_k$. Finally, the equilibrium cost of agent $i \in \mathcal{I}=\{1, \ldots, I \}, $  is 
\begin{equation*}
    f_{i0}(x_0, \bar{x}_0) = {\alpha}_{i0} {\mathbb{E}[(x_0-\bar{x}_0)^{2o}] }+ \bar{\alpha}_{i0} \bar{x}_0^{2p},
\end{equation*}
completing the proof. \qed
\end{proof}

\section{Numerical Examples}
\label{sec:examples}

In this Section we show some numerical examples that allow illustrating the game theoretical results presented in this paper.

\begin{figure}[t!]
    \centering
    \begin{tabular}{cc}
       \includegraphics[width=0.5\columnwidth]{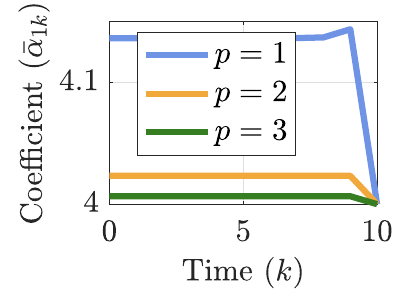}  &
       \hspace{-0.5cm}
       \includegraphics[width=0.5\columnwidth]{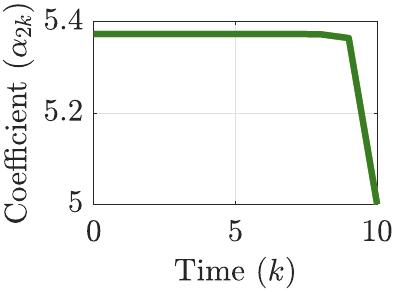}  \\
       \includegraphics[width=0.5\columnwidth]{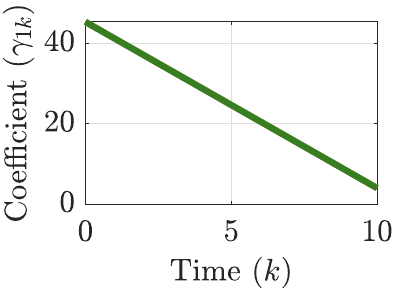}  &
       \hspace{-0.5cm}
       \includegraphics[width=0.5\columnwidth]{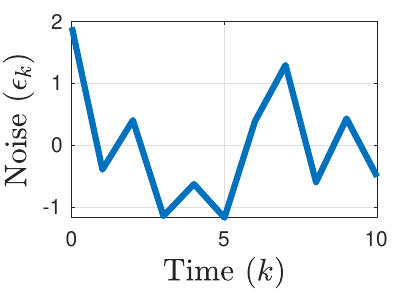}  \\
       \includegraphics[width=0.5\columnwidth]{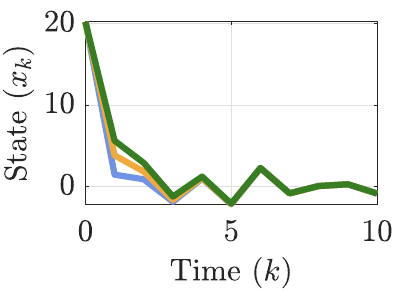}  &
       \hspace{-0.5cm}
       \includegraphics[width=0.5\columnwidth]{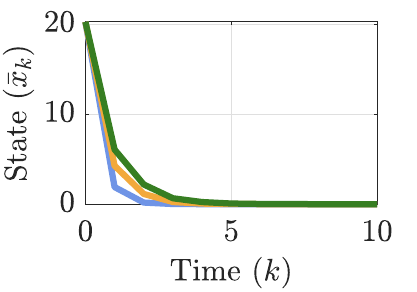}  \\
         \includegraphics[width=0.5\columnwidth]{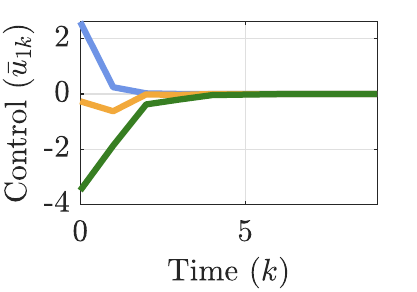}  &
       \hspace{-0.5cm}
       \includegraphics[width=0.5\columnwidth]{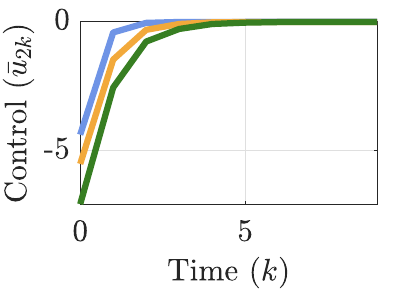}  \\
       \includegraphics[width=0.5\columnwidth]{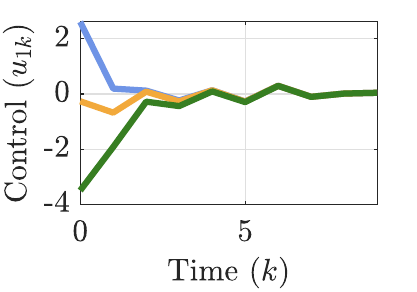}  &
       \hspace{-0.5cm}
       \includegraphics[width=0.5\columnwidth]{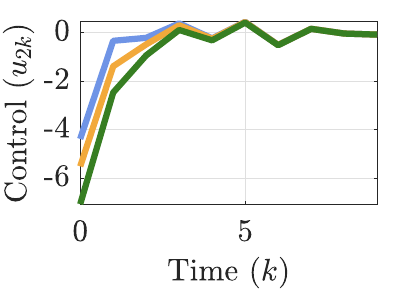}  
    \end{tabular}
    \caption{Results example corresponding to Proposition \ref{propos:proposition_10} and Proposition \ref{propos:proposition_11}.}
    \label{fig:p3_4_games}

    \vspace{-0.3cm}
\end{figure}

\vspace{-0.3cm}
\subsection{Example 1: Higher Order Costs}

This example corresponds to the results presented in Proposition \ref{propos:proposition_8} and Proposition \ref{propos:proposition_9} considering two players. We study the effect of $p$ on $\bar{\alpha}_{1k}$ and $\bar{\alpha}_{2k}$, $x$, and strategies $\bar{u}_{1k}$ and $\bar{u}_{2k}$. Let us consider a system with $\bar{a}_k = 1$, $\bar{b}_{1k} = -2$, $\bar{b}_{2k} = 2$, $\bar{q}_{1k} = \bar{q}_{1N} = 4$, $\bar{q}_{2k} = \bar{q}_{2N} = 5$, $\bar{r}_{1k} = 6$, $\bar{r}_{2k} = 7$ and $\bar{x}_0 = 10$; and the terminal time $N=7$. The results are presented in Figure \ref{fig:p1_2_games}. It can be seen the backward recursive equation evolution $\bar{\alpha}_{1k}$ depending on the parameter $p$. For all the cases, and as expected, we observe that the state goes to zero minimizing the cost. The equilibrium strategies for both players $\bar{u}_{1k}^*$ and $\bar{u}_{2k}^*$ are presented in Figure \ref{fig:p1_2_games}. Note that they exhibit a completely different behavior as their cost functions have different parameters, and in addition, the system parameters $\bar{b}_{1k}$ and $\bar{b}_{2k}$ have different sign.

\vspace{-0.3cm}
\subsection{Example 2: Variance-Aware Higher-Order Costs}

This example corresponds to the results presented in Proposition \ref{propos:proposition_10} and Proposition \ref{propos:proposition_11} with two players. We study the effect of $p$ on $\bar{\alpha}_{ik}$, ${\alpha}_{ik}$, $\gamma_{ik}$, $x$, and $\bar{u}_{ik}$, for all $i \in \{1,2\}$. Let us consider a system with $\bar{a}_k = 1$, $\bar{b}_{1k} = -2$, $\bar{b}_{2k} = 3$, $\epsilon$ as in Figure \ref{fig:p3_4_games} with unitary variance $\sigma=1$, and $\bar{q}_{1k} = \bar{q}_{1N} = 4$, $\bar{q}_{2k} = \bar{q}_{2N} = 5$, $\bar{r}_{1k} = 6$, $\bar{r}_{2k} = 7$, ${q}_{1k} = {q}_{1N} = 4$, ${q}_{2k} = {q}_{2N} = 5$, ${r}_{1k} = 6$, ${r}_{2k} = 7$ and $\bar{x}_0 = x_0 = 20.25$; and the terminal time $N=10$. The results are presented in Figure \ref{fig:p3_4_games}. We can observe that the recursive equations for $\alpha_{2k}$ and $\bar{\gamma}_{1k}$ are the same for all the different $p$ values as shown in Proposition \ref{propos:proposition_11}. Note that the optimal strategy for the first player dramatically changes the trend as we change the $p$ value.

\vspace{-0.3cm}
\subsection{Example 3: Risk-Aware Games under Multiplicative Noise}

\begin{figure}[t!]
    \centering
    \begin{tabular}{cc}
       \includegraphics[width=0.5\columnwidth]{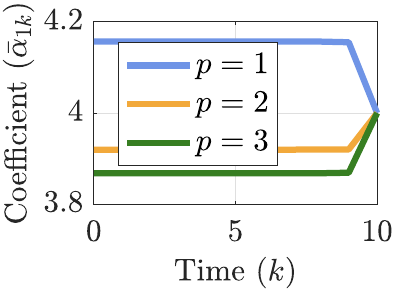}  &
       \hspace{-0.5cm}
       \includegraphics[width=0.5\columnwidth]{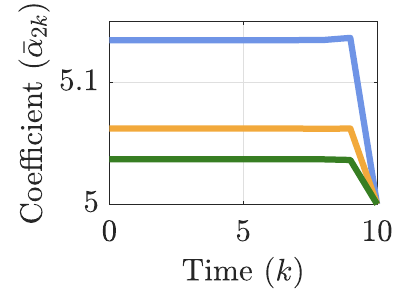}  \\
       \includegraphics[width=0.5\columnwidth]{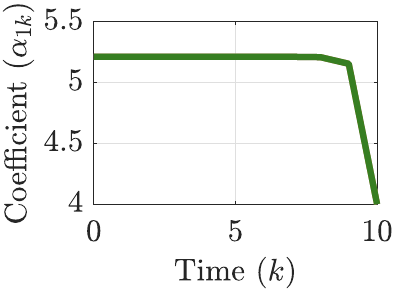}  &
       \hspace{-0.5cm}
       \includegraphics[width=0.5\columnwidth]{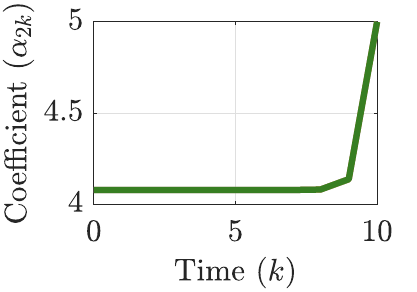}  \\
       \includegraphics[width=0.5\columnwidth]{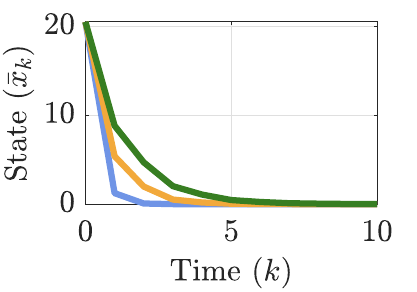}  &
       \hspace{-0.5cm}
       \includegraphics[width=0.5\columnwidth]{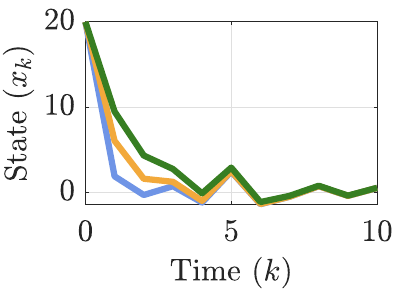}  \\
         \includegraphics[width=0.5\columnwidth]{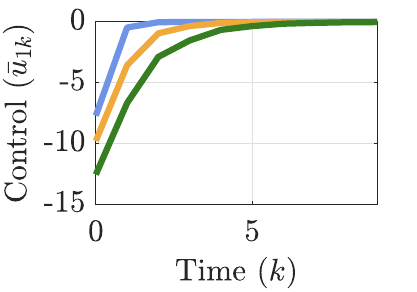}  &
       \hspace{-0.5cm}
       \includegraphics[width=0.5\columnwidth]{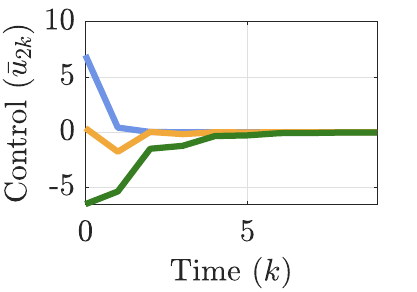}  \\
       \includegraphics[width=0.5\columnwidth]{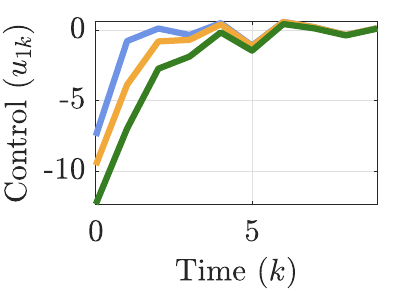}  &
       \hspace{-0.5cm}
       \includegraphics[width=0.5\columnwidth]{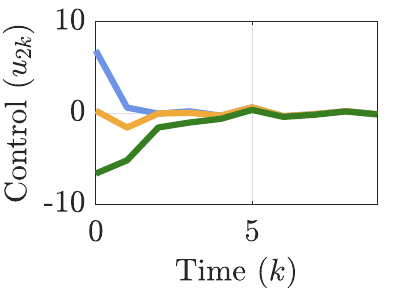}  
    \end{tabular}
    \caption{Results example corresponding to Proposition \ref{propos:proposition_12} and Proposition \ref{propos:proposition_13}.}
    \label{fig:p5_6_games}

    \vspace{-0.3cm}
\end{figure}

This example corresponds to the results presented in Proposition \ref{propos:proposition_12} and Proposition \ref{propos:proposition_13}. We study the effect of $p$ on $\bar{\alpha}_{ik}$, ${\alpha}_{ik}$, $x$, and $\bar{u}_{ik}$, for all $i \in \{1,2\}$. Let us consider a system with $\bar{a}_k = 1$, $\bar{b}_{1k} = 1.5$, $\bar{b}_{2k} = -1.1$, $\epsilon$ as in Figure \ref{fig:p3_4_games} with unitary variance $\sigma=1$, and $\bar{q}_{1k} = \bar{q}_{1N} = 4$, $\bar{q}_{2k} = \bar{q}_{2N} = 5$, $\bar{r}_{1k} = \bar{r}_{2k} = 1$, ${q}_{1k} = {q}_{1N} = 5$, ${q}_{2k} = {q}_{2N} = 4$, ${r}_{1k} ={r}_{2k} = 1$ and $\bar{x}_0 =20.5$, $x_0 =20$; and the terminal time $N=10$.
Figure \ref{fig:p5_6_games} presents the results for different $p$ values. We observe that the system state is driven to zero as desired according to the cost functional. In addition, this controller is minimizing higher-order terms as we increase the $p$ value. Interestingly, we observe a smaller variations in the optimal strategies as we increase the $p$ value. For the strategy of the second player, we observe a different optimal behavior depending on the $p$ value.

\section{Concluding remarks}
\label{sec:conclusions}

We  investigated higher-order costs in the context of discrete time mean-field-type  game theory, addressing limitations in traditional quadratic cost models and offering solutions for various scenarios.  We presented semi-explicit solutions for different problem types, including  selfish agents with higher-order costs, selfish agents with variance-aware higher-order costs, and risk-aware games under additive and multiplicative noise.  These solutions provide a valuable tool for understanding and controlling multi-agent systems that involve non-linear penalties and risk-sensitive behavior.   

A potential avenue for future work involves exploring the application of these semi-explicit solutions to specific real-world problems, such as designing efficient energy systems, developing robust wireless networks, or optimizing interactive blockchain technologies.  Additionally, extending the analysis to include other risk measures, such as those incorporating skewness and kurtosis, expectile, extremile, could further enhance the applicability of the framework.

\bibliographystyle{unsrt}
\bibliography{references.bib}

\begin{thebibliography}{10}

\bibitem{Coopetitive}
J.~Barreiro-Gomez, T.~E. Duncan, and H.~Tembine.
\newblock Co-opetitive linear-quadratic mean-field-type games.
\newblock {\em IEEE Transactions on Cybernetics}, 50(12):5089 -- 5098, 2020.

\bibitem{berge}
N.~Toumi, J.~Barreiro-Gomez, T.~E. Duncan, and H.~Tembine.
\newblock Berge equilibrium in linear-quadratic mean-field-type games.
\newblock {\em Journal of the Franklin Institute}, 357(15):10861--10885, 2020.

\bibitem{Jovanovic_2}
B.~Jovanovic and R.~W. Rosenthal.
\newblock Anonymous sequential games.
\newblock {\em Journal of Mathematical Economics}, 17(1):77--87, 1988.

\bibitem{Lions}
J.~M. Lasry and P.~L. Lions.
\newblock Mean field games.
\newblock {\em Japanese Journal of Mathematics}, 2(2007):229--260, 2007.

\bibitem{Huang}
M.~Huang, R.~P. Malham\'e, and P.~E. Caines.
\newblock Large population stochastic dynamic games: closed-loop mckean-vlasov
  systems and the nash certainty equivalence principle.
\newblock {\em Communications in information and systems}, 6(2006):221--251,
  2006.

\bibitem{MAS_HT}
H.~Tembine.
\newblock Master adjoint systems in mean-field-type games.
\newblock {\em Communications in Information and Systems}, 21(4):623--650,
  2021.

\bibitem{nonlinearMFTG}
J.~Barreiro-Gomez, T.~E. Duncan, B.~Pasik-Duncan, and H.~Tembine.
\newblock Semi-explicit solutions to some non-linear non-quadratic
  mean-field-type games: A direct method.
\newblock {\em IEEE Transactions on Automatic Control}, 65(6):2582--2597, 2020.

\bibitem{refmain0}
H.~Tembine.
\newblock Risk-sensitive mean-field-type games with {Lp}-norm drifts.
\newblock {\em Automatica}, 59(2015):224--237, 2015.

\bibitem{refmain1}
B.~Djehiche, A.~Tcheukam, and H.~Tembine.
\newblock Mean-field-type games in engineering.
\newblock {\em AIMS Electronics and Electrical Engineering}, 1(1):18--73, 2017.

\bibitem{refmain2}
H.~Tembine.
\newblock Mean-field-type games.
\newblock {\em AIMS Math}, 2(4):706--735, 2017.

\bibitem{Higher-order-control}
J.~Barreiro-Gomez, T.~E. Duncan, B.~Pasik-Duncan, and H.~Tembine.
\newblock Semi-explicit solution of some discrete-time mean-field-type control
  with higher-order costs.
\newblock {\em arXiv preprint:2505.04112}, 2025.

\bibitem{ref001}
H.M. Markowitz.
\newblock Portfolio selection.
\newblock {\em The Journal of Finance}, 7(1):77--91, 1952.

\bibitem{ref002}
C.F. Roos.
\newblock A mathematical theory of competition.
\newblock {\em American Journal of Mathematics}, 47(3):163--175, 1925.

\bibitem{ref003}
C.F. Roos.
\newblock A dynamic theory of economics.
\newblock {\em Journal of Political Economy}, 35:632--656, 1927.

\bibitem{ref004}
H.M. Markowitz.
\newblock The utility of wealth.
\newblock {\em Journal of Political Economy}, 60:151--158, 1952.

\bibitem{ref005}
H.M. Markowitz.
\newblock {\em Portfolio Selection: Efficient Diversification of Investments}.
\newblock John Wiley \& Sons. New York, NY, USA, 1959.

\bibitem{ref006}
A.~Bensoussan, J.~Frehse, and S.C.P. Yam.
\newblock {\em Mean Field Games and Mean Field Type Control Theory}.
\newblock Springer: Berlin, Germany, 2013.

\bibitem{ref04}
A.~Tcheukam Siwe and H.~Tembine.
\newblock On the distributed mean-variance paradigm.
\newblock In {\em 13th International Multi-Conference on Systems, Signals \&
  Devices (SSD)}, pages 607--612, Leipzig, Germany, 2016.
\newblock DOI: 10.1109/SSD.2016.7473660.

\bibitem{ref01}
A.~K. Cisse and H.~Tembine.
\newblock Cooperative mean-field type games.
\newblock {\em IFAC Proceedings}, 47(3):8995--9000, 2014.

\bibitem{ref02}
H.~Tembine.
\newblock Uncertainty quantification in mean-field-type teams and games.
\newblock In {\em 54th IEEE Conference on Decision and Control (CDC)}, Osaka,
  Japan, 2015.
\newblock DOI: 10.1109/CDC.2015.7402909.

\bibitem{ref03}
A.~Tcheukam Siwe and H.~Tembine.
\newblock Network security as public good: A mean-field-type game theory
  approach.
\newblock In {\em 13th International Multi-Conference on Systems, Signals \&
  Devices (SSD)}, pages 601--606, Leipzig, Germany, 2016.
\newblock DOI: 10.1109/SSD.2016.7473659.

\bibitem{ref0}
T.~E. Duncan and H.~Tembine.
\newblock Linear–quadratic mean-field-type games: A direct method.
\newblock {\em Games}, 9(1), 2018.

\bibitem{ref0b}
J.~Barreiro-Gomez and H.~Tembine.
\newblock Blockchain token economics: A mean-field-type game perspective.
\newblock {\em IEEE Access}, 7(2019):64603--64613, 2019.

\bibitem{ref1}
J.~Barreiro-Gomez and H.~Tembine.
\newblock {\em Mean-Field-Type Games for Engineers}.
\newblock CRC Press, Inc., 2021.

\bibitem{Toolbox}
J.~Barreiro-Gomez and H.~Tembine.
\newblock A {MatLab}-based mean-field-type games toolbox: Continuous-time
  version.
\newblock {\em IEEE Access}, 7:126500--126514, 2019.

\bibitem{ref2}
J.~Barreiro-Gomez, T.~E. Duncan, and H.~Tembine.
\newblock Discrete-time linear-quadratic mean-field-type repeated games:
  Perfect, incomplete, and imperfect information.
\newblock {\em Automatica}, 112(2020):108647, 2020.

\bibitem{ref3}
J.~Barreiro-Gomez, T.~E. Duncan, and H.~Tembine.
\newblock Linear-quadratic mean-field-type games-based stochastic model
  predictive control: A microgrid energy storage application.
\newblock In {\em American Control Conference (ACC)}, pages 3224--3229,
  Philadelphia, PA, USA, 2019.

\bibitem{waterref}
J.~Barreiro-Gomez and H.~Tembine.
\newblock Mean-field-type model predictive control: An application to water
  distribution networks.
\newblock {\em IEEE Access}, 7:135332--135339, 2019.

\bibitem{ref4}
R.~Liu, Y.~Li, and X.~Liu.
\newblock Linear-quadratic optimal control for unknown mean-field stochastic
  discrete-time system via adaptive dynamic programming approach.
\newblock {\em Neurocomputing}, 282:16--24, 2018.

\bibitem{ref4b}
B.~Djehiche, J.~Barreiro-Gomez, and H.~Tembine.
\newblock Price dynamics for electricity in smart grid via mean-field-type
  games.
\newblock {\em Dynamic Games and Applications}, 10(2020):798--818, 2020.

\bibitem{ref4c}
J.~Barreiro-Gomez, T.~E. Duncan, and H.~Tembine.
\newblock Linear-quadratic mean-field-type games with multiple input
  constraints.
\newblock {\em IEEE Control Systems Letters}, 3(3):511--516, 2019.

\bibitem{viswa_2023}
P.~S. Mohapatra and P.~V. Reddy.
\newblock Linear-quadratic mean-field-type difference games with coupled affine
  inequality constraints.
\newblock {\em IEEE Control Systems Letters}, 7:1987 -- 1992, 2023.

\bibitem{ref5b}
J.~Barreiro-Gomez, T.~Duncan, B.~Pasik-Duncan, and H.~Tembine.
\newblock Semi-explicit solutions to some nonlinear nonquadratic
  mean-field-type games: A direct method.
\newblock {\em IEEE Transactions on Automatic Control}, 65(6):2582--2597, 2020.

\bibitem{ref5c}
J.~Barreiro-Gomez, B.~Djehiche, T.~E. Duncan, B.~Pasik-Duncan, and H.~Tembine.
\newblock Fractional mean-field-type games under non-quadratic costs: A direct
  method.
\newblock In {\em IEEE 58th Conference on Decision and Control (CDC)}, pages
  293--298, Nice, France, 2019.

\bibitem{ref5}
Z.~El~Oula Frihi, J.~Barreiro-Gomez, S.~Eddine Choutri, and H.~Tembine.
\newblock Hierarchical structures and leadership design in mean-field-type
  games with polynomial cost.
\newblock {\em Games}, 11(30), 2020.

\bibitem{ref7}
H.~Tembine.
\newblock Distributed planning in mean-field-type games.
\newblock In {\em 21st IFAC World Congress}, pages 2183--2188, Berlin, Germany,
  2020.

\bibitem{ref6}
T.~E. Duncan, B.~Pasik-Duncan, and H.~Tembine.
\newblock Mean-field-type games driven by rosenblatt processes.
\newblock In {\em 10th International Conference on Control, Decision and
  Information Technologies (CoDIT)}, pages 982--987, Vallette, Malta, 2024.

\end{thebibliography}

\end{document}